\documentclass{article}
\usepackage[utf8]{inputenc}
\usepackage[left=2cm, top=2cm, right=2cm, bottom=2cm]{geometry}
\usepackage{amsmath,amssymb,amsfonts}
\usepackage{url}
\usepackage[caption=false,font=footnotesize]{subfig}
\usepackage{algcompatible}
\usepackage{algorithm}
\usepackage{cite}
\usepackage[pdftex]{graphicx}

\usepackage[hidelinks]{hyperref}
\usepackage{bm}
\usepackage{color}
\usepackage{dsfont}
\usepackage{wasysym}
\usepackage{array}
\usepackage{multirow}
\usepackage{fancyhdr,lastpage}
\providecommand{\keywords}[1]{{\textbf{Key words---}} #1}

\usepackage{amsthm}
\theoremstyle{plain}
\newtheorem{theorem}{Theorem}
\newtheorem{lemma}{Lemma}

\newtheorem{proposition}{Proposition}

\DeclareMathOperator*{\argmin}{arg\,min}

\DeclareMathOperator{\E}{\mathbb{E}}

\DeclareMathOperator{\card}{card}

\def\Halmos{\mbox{\quad$\square$}}

\newcommand{\RR}{\mathbb{R}}

\newcommand{\UU}{\mathcal{U}}

\newcommand{\XX}{\mathbf{X}}
\newcommand{\xx}{\mathbf{x}}
\newcommand{\xs}{\bm{\xi}}
\newcommand{\uu}{\mathbf{u}}

\usepackage{natbib}
 \bibpunct[, ]{(}{)}{,}{a}{}{,}%
 \def\bibfont{\small}%
 \def\BIBand{and}%

\title{Optimally Scheduling Public Safety Power Shutoffs}

\author{Antoine Lesage-Landry, Félix Pellerin, Joshua A. Taylor, and Duncan S. Callaway\footnote{A. Lesage-Landry and F. Pellerin are with the Department of Electrical Engineering, Polytechnique Montréal \& GERAD, Montréal, QC, Canada, H3T 1J4. e-mail: \{\texttt{antoine.lesage-landry}, \texttt{felix.pellerin}\}\texttt{@polymtl.ca}. J.A. Taylor is with The Edward S. Rogers~Sr. Department of Electrical \& Computer Engineering, University of Toronto, Toronto, Ontario, Canada, M5S 3G4. e-mail: \texttt{josh.taylor@utoronto.ca}. D.S. Callaway is with the Energy \& Resources Group, University of California, Berkeley, CA, USA 94720. e-mail: \texttt{dcal@berkeley.edu}.}}

\date{}

\begin{document}
\pagestyle{fancy}

\maketitle

\begin{abstract}
In an effort to reduce power system-caused wildfires, utilities carry out public safety power shutoffs (PSPS) in which portions of the grid are de-energized to mitigate the risk of ignition. The decision to call a PSPS must balance reducing ignition risks and the negative impact of service interruptions. In this work, we consider three PSPS scheduling scenarios, which we model as dynamic programs. 
In the first two scenarios, we assume that $N$ PSPSs are budgeted as part of the investment strategy. In the first scenario, a penalty is incurred for each PSPS declared past the $N^\text{th}$ event. In the second, we assume that some costs can be recovered if the number of PSPSs is below $N$ while still being subject to a penalty if above $N$. In the third, the system operator wants to minimize the number of PSPS such that the total expected cost is below a threshold. We provide optimal or asymptotically optimal policies for each case, the first two of which have closed-form expressions. Lastly, we establish the applicability of the first PSPS model's policy to critical-peak pricing, and obtain an optimal scheduling policy to reduce the peak demand based on weather observations. 
\end{abstract}

\bigskip

\section{Introduction}
Electric power systems have caused a number of recent wildfires, for example, in California, USA (2017, 2018), in Texas, USA (2011), and in Victoria, AU (2009)~\citep{WRR,teague2010final,russell2012distribution,abatzoglou2020population,jazebi2019review,rhodes2020balancing}.
The consequences of these events are exacerbated by the extreme weather conditions during which they often occur~\citep{miller2017electrically,keeley2018historical}.
To mitigate fires, utilities like Northern California's Pacific Gas \& Electricity (PG\&E) have implemented precautionary grid de-energization events, i.e., intentional load shedding~\citep{PSPSpolicies,abatzoglou2020population,rhodes2020balancing}. These events are referred to as public safety power shutoffs (PSPSs) by PG\&E and in the literature. However, de-energizing part of the grid has important consequences on the affected areas, e.g., temporary business and school closures, medical baseline customers, and loss of revenue~\citep{abatzoglou2020population,rhodes2020balancing,PSPSmedicalbaseline}. 

In this work, we formulate dynamic programs for PSPS scheduling {in a specific geographic area of a power grid}. The models use observations of natural phenomena to balance estimated wildfire risks with the cost of de-energization.
Our models assume that initial mitigation investments are made, for example, tree-trimming, better line insulation, and under-grounding of lines~\citep{PGEWildfireMit}. We study the following three scenarios:
\begin{enumerate}
    \item $N$ PSPSs are budgeted. We assume that the system operator is allowed to exceed the budget but is penalized for each extra PSPS{, e.g., due to contractual agreements regarding equipment rental and associated extension costs, and to fixed on-the-ground staff expenditures (Section~\ref{sec:add_PSPS}).}
    \item The system operator can recover part of their incurred costs if the number of PSPS is below the budget. {For example, community resource centers which are used to support the population during outages do not need to be deployed and costs of personnel salary and food are avoided, and portable back-up generators or hotel ticket vouchers do not need to be provided to medical baseline consumer~\citep{PSPSfactsheer} (Section~\ref{sec:adjustment}).}
    \item The system operator minimizes the number of PSPSs such that the total expected cost is below a threshold that depends on {the wildfire mitigation strategy investment. This scenario applies to, for example, jurisdictions where the system operator is not able to forecast its needs in term of PSPSs, e.g., due to limited experience in regions having not faced powerline-caused wildfire in the recent past (Section~\ref{sec:minPSPS}).}
\end{enumerate}
For each case, we provide the optimal scheduling policy or an asymptotically equivalent form of it.

We show that first model's asymptotic policy is optimal for critical-peak pricing (CPP)~\citep{siano2014demand}.
Loads enrolled in CPP see very high electricity prices when the utility needs to reduce demand, e.g., due to extreme hot or cold weather, and pay a discounted rate at other times~\citep{herter2007residential,chen2013optimal}. 
CPP events are called a day ahead by the utility~\citep{vardakas2014survey}. 
The total number of CPP events must not exceed a contracted limit, which we show is
equivalent to the first PSPS scenario's relaxed model. We use its optimal policy for scheduling CPP events based on natural phenomenon observations.

\paragraph*{Related work.}
Scheduling PSPSs is a topic of relatively recent interest. For example, PG\&E only started using large-scale PSPSs in 2019~\citep{abatzoglou2020population}. 
{In their work,~\cite{rhodes2020balancing} formulated an optimal PSPS scheduling problem in which the wildfire risks are balanced with the costs of power outages. In~\citep{astudillo2022managing}, a multiperiod optimal PSPS scheduling problem was formulated as a deterministic mixed-integer linear program. Energy storage was considered in their formulation to increase grid reliability when subject to PSPSs. In~\citep{rhodes2022co}, the authors used a rolling horizon formulation to include new information about wildfires. The cost of line re-energization was also considered in~\citep{rhodes2022co}. Fairness considerations were added in~\citep{kody2022sharing}'s PSPS scheduling formulation which is expressed as mixed-integer linear model predictive control problem. Learning-based approaches were used in~\citep{umunnakwe2022data} to estimate wildfire risks, the outputs of which were then used to plan preventive line outages or load shedding events. Reference~\citep{hong2022data} developed a data-driven decision-making method for efficiently planning load shedding and PSPS events based on a data set generated from optimal power flow computations. In related work,~\citep{kody2022optimizing} formulated an optimal infrastructure investment strategy to reduce both wildfire risks and PSPS impacts.} Our work differs from the above in that the decision-making is sequential and natural phenomena are modeled as stochastic processes.

We now review the literature on CPP scheduling. Reference~\citep{joo2007option} combines market price prediction and swing options to call CPPs in a way that maximizes the profit of the system operator. In~\citep{zhang2009optimal}, load price elasticity is used to balance utility profits and load costs when calling CPP. CPP scheduling was framed as a dynamic program in~\citep{chen2013optimal} and threshold policies for different types of CPPs were proposed. Our work differs from~\citep{chen2013optimal} in that multiple natural phenomena can be considered in the decision-making process. 
We also prove that the threshold policy is optimal for fixed-period CPP, which was not done in~\citep{chen2013optimal}. However, we remark that~\citep{chen2013optimal} can accommodate a number of variants that our approach cannot, e.g., variable peak pricing, variable-period CPP, and multiple-group CPP.

We obtain our optimal policies by adapting dynamic programming-based techniques for the multiple secretary problem~\citep{kleinberg2005multiple,babaioff2007knapsack,arlotto2019uniformly}. This problem is an extension of the secretary problem first published in a scientific journal by~\citep{lindley1961dynamic} in which the $N \geq 2$ best candidates with random abilities must be selected. References~\citep{freeman1983secretary,ferguson1989solved} provide detailed surveys of the secretary problem.

The remainder of the paper is as follows. The PSPS scheduling problem is introduced in Section~\ref{sec:psps_planning}. Our three models and their respective policies are presented in Sections~\ref{sec:add_PSPS},~\ref{sec:adjustment}, and~\ref{sec:minPSPS}. We provide detailed and original derivations of Scenario 1 and 2's policies in Sections~\ref{sec:analysis_costpenalty} and~\ref{ssec:analysis_adj}, and leverage the results of~\citep{chen2004dynamic,chen2007non} to provide an optimal policy for Scenario 3 in Section~\ref{sec:analysis_scenario3}. In Section~\ref{sec:cpp}, we apply the results of Section~\ref{sec:add_PSPS} to CPP. Section~\ref{sec:num} presents numerical simulations for the first two PSPS and the CPP models. We conclude in Section~\ref{sec:conclu}.

\section{PSPS Scenarios}
\label{sec:psps_planning}
We now present our novel PSPS scheduling models and their respective optimal or asymptotically optimal policies.
We consider a sequence of days indexed by $t = 1, 2, \ldots, T$. Let $C_I > 0$ be the initial investment in infrastructure upgrades made to reduce wildfire risks. We consider a single geographic area prone to powerline-caused wildfires, e.g., supplied by aging or hardly accessible transmission lines.
Let $u_{t} = 1$ be the decision taken at day $t$ to call for a PSPS in this region for the next day, i.e., at $t+1$. Let $u_{t} = 0$ be the complementary decision to operate the grid normally for the next day. Let $\mathbf{u} = \left(u_1, u_2, \ldots, u_{T} \right)^\top \in \left\{0, 1 \right\}^T$. For the first two scenarios, we let $0 \leq k < N$ be the remaining PSPS budget at round $t$:
\[
k = N - \sum_{i=0}^{t-1} u_i.
\]

Let $\XX_t \in \RR^n$ be a random vector made of observations and measurements from $n \in \mathbb{N}$ different natural phenomena and other factors impacting wildfire ignitions. The variable $\XX_t$ is observed at the end of day $t$. The entries of $\XX_t$ can represent factors such ambient temperature or wind speed~\citep{PSPSfactsheer}. 
We assume that $\XX_t$ is a discrete-time, finite-state Markov process with known transition probabilities. Let $\mathcal{X} \subset \RR^n$ be the state space of $\XX_t$. Let $\mathbf{P}\in \RR^{n \times n}$ be the transition matrix. We denote the probability of moving from state $\mathbf{y} \in \mathcal{X}$ to $\mathbf{x} \in \mathcal{X}$ by $\Pr\left[\left.\XX_{t+1} = \xx \right| \XX_{t} = \mathbf{y}  \right]$.

For $t=1,2,\ldots, T$, let $f_t: \RR^n \mapsto \left\{0,1\right\}$ be a function mapping $\XX_t$ to a binary output indicating the high risk ($f_t(\XX_t) = 1$) or low risk ($f_t(\XX_t) = 0$) of ignition given that the lines are energized during day $t$. The function $f_t$ is assumed to be deterministic, known, and set according to the system operator needs and knowledge of its geographical area and infrastructure. For example, $f_t$ could take the form:
\begin{equation}
f_t\left( \XX_t \right) = \begin{cases}
1, &\text{ if } \XX_t(i) \geq \bm{\Delta}_t(i),  \ i=1,2,\ldots,n \\
0, &\text{ otherwise, }
\end{cases}\label{eq:risk_thres}
\end{equation}
where $\bm{\Delta}_t(i)$ is the risk threshold for observation $i$. This function is inspired by PG\&E's `minimum fire conditions'~\citep{PSPSpolicies}. For example, PG\&E uses thresholds like: $40$ km/h for sustained wind speed, $72$ km/h for wind gust, and an air humidity below $20\%$~\citep{PSPSfactsheer}. Hence, if the entries of $\XX_t$ are greater or equal to these values for the wind speed and gust, and less or equal to for the humidity, then $f_t(\XX_t) = 1$.

Let $A_t >0$ be the estimated cost of powerline-caused wildfires at time $t$. This represents loss of life, material damage, reconstruction and repair, and service interruption. The powerline-induced wildfire damage cost at $t$ is then given by $A_t f_t(\mathbf{X}_t)$, i.e., the estimated cost weighted by the risk indicator function.
Let $a_t >0$ be the revenue loss of the utility when de-energizing the grid at time $t$.  
Let $s_{1,t} > 0$ and $s_{2,t} >0$ be, respectively, the cost of de-energizing and re-energizing the grid. The scalar coefficient $s_{1,t}$ includes, for example, the costs of crews and special equipment dispatched to de-energize the grid. The scalar coefficient $s_{2,t}$ also includes the costs related to the personnel on the field in addition to the cost of specialized equipment required to inspect the network before re-energization~\citep{PSPSpolicies}. We model the cost of a de/re-energizing as: $\max\left\{s_{1,t+1} \left( u_{t} - u_{t-1} \right), s_{2,t+1}\left(u_{t-1} - u_{t} \right)  \right\}$. We assume that the grid is energized before and after the time range, i.e., $u_{0} = u_{T+1} = 0$.

\subsection{Scenario~1: Cost Penalty for Additional PSPSs}
\label{sec:add_PSPS}
We assume that the wildfire mitigation strategy does not strictly rely on infrastructure upgrades, and PSPSs can be scheduled to reduce the strategy's total costs, $C_I$. In the first two scenarios, we assume that $C_I$ is decomposed into two components: (i) a reduced budget for infrastructure upgrades and PSPS-induced financial losses when $N$ PSPSs have been planned for, which we denote $\overline{C}_I(N) < C_I$, and (ii) the operational and/or damage costs associated with calling a PSPS or not, at each day. We can use, for example, the value of lost load~\citep{kariuki1996evaluation,willis1997electricity,ratha2013value} for a day to model the financial losses from the former component.

In Scenario~1, we assume that additional costs are incurred for each PSPS after the budget, $N$, has been reached. These costs represent, for example, the value of lost load per extra day without power and the cost of personnel and equipment contract extensions. These costs are not recovered if the number of PSPS is below $N$ because of contracted obligations, e.g., surveillance personnel and hired equipment like helicopters for post-PSPS inspections~\citep{PGEWildfireMit}. Let $\gamma > 0$ be the cost penalty per extra power shutoff. 
We model the additional costs as: $\gamma \max\left\{ 0, \sum_{t=1}^T u_{t} - N \right\}$.

The total cost with penalty for Scenario~1 is:
\begin{equation}
\begin{aligned}
c_\text{p}\left(N, \uu \right) &= \sum_{t=1}^T \left( a_{t+1} u_{t} +  A_{t+1} f_{t+1}\left(\XX_{t+1}\right) (1-u_{t}) + \max\left\{{s_{1,t+1} \left( u_{t} - u_{t-1} \right), s_{2,t+1}\left(u_{t-1} - u_{t} \right) } \right\}\right)\\
&\quad + \gamma\max \left\{0, \sum_{t=1}^T u_{t} - N \right\} + \overline{C}_I(N).
\end{aligned}
\label{eq:total_cost}
\end{equation}
Our objective is to minimize the expectation of $c_\text{p}\left(N, \uu \right)$ given $\overline{C}_I(N)$. We hereon omit $\overline{C}_I(N)$ in both Scenario 1's and 2's objective functions because it does not influence their optima. 
The PSPS scheduling problem is expressed as the following problem, which we later rewrite as a dynamic program:
\begin{equation}
\begin{aligned}
& \min_{\substack{u_t\\ t=1,2,\ldots,T}} & & \E \left[ c_\text{p}\left(N, \uu \right) \right]\\
& \text{subject to}
& & u_t \in \left\{0, 1 \right\}\\
& & & u_{0} = u_{T+1} = 0,
\end{aligned}
\label{eq:min_unconstrained}
\end{equation}
Note that the expectation is taken with respect to the random variables $\XX_t$, $t=1,2,\ldots, T$. 

Problem~\eqref{eq:min_unconstrained} is approximately solved by the following closed-form, easily-implementable policy. Proposition~\ref{prop:policy_psps}'s policy is optimal for an asymptotically equivalent form of Problem~\eqref{eq:min_unconstrained}. This will be shown in Section~\ref{sec:relax}.

\begin{proposition}
Consider problem~\eqref{eq:min_unconstrained}. 
Let $d=T-t$ and the current day's observations be $\XX_{d} = \xx$. Then, $u_d =1$ if
\begin{align*}
\E\left[\left.f_{d-1}\left(\XX_{d-1} \right)\right| \XX_{d} = \xx \right] \geq& \frac{1}{A_{d-1}} \big(g_{d-1}(k-1|1,\xx) - g_{d-1}(k|0,\xx) + a_{d-1}\\
&\qquad\qquad + \left(1-u_{d+1}\right) s_{{1},d-1} - u_{d+1} s_{{2},d-1} \big),
\end{align*}
where
\begin{equation}
\begin{aligned}
g_{d}(k|u_{d+1} = u,\XX_{d+1} =\xx) &= \sum_{\xs \in \mathcal{X}} \min \left\{ g_{d-1}(k|0,\xs) +  \E\left[ \left.A_{d-1}f_{d-1}\left(\XX_{d-1}\right)\right| \XX_d = \xs \right] + u s_{{2},d-1} ,\right.\\
&\qquad\qquad\quad \left. g_{d-1}(k-1|1,\xs) + a_{d-1} + (1-u) s_{{1},d-1}\right\}  \Pr\left[\left.\XX_d = \xs \right| \XX_{d+1} = \xx \right], \label{eq:def_g}
\end{aligned}
\end{equation}
with
\begin{align}
g_{0}(k|u,\xx) &= u s_{{2},{T+1}} + \sum_{\xs \in {\mathcal{X}}}\E\left[ \left.A_{{T+1}}f_{{T+1}}\left(\XX_{{T+1}}\right)\right| \XX_{{T}} = \xs \right]\Pr\left[\left.\XX_{T} = \xs \right| \XX_{{T-1}} = \xx \right] \label{g_0d}\\
g_{d}(0|u,\xx) &=u s_{{2},d-1} + \sum_{i=0}^d\sum_{\xs \in \mathcal{X}}\E\left[ \left.A_{T-i+1}f_{T-i+1}\left(\XX_{T-i+1}\right)\right| \XX_{T-i} = \xs \right]\mathbf{P}^{(i+1)}_{\xx,\xs} ,
\label{eq:g_0}
\end{align}
for all $\xx \in \mathcal{X}, u \in \left\{0,1 \right\}, k \geq 1 \text{ and } d \geq 1$, and where $\mathbf{P}^{(i+1)}_{\xx,\xs} = \Pr\left[\left. \XX_{n+i+1} = \xs \right| \XX_n = \xx \right]$ for any $n$.
\label{prop:policy_psps}
\end{proposition}

\subsection{Scenario~2: Cost Adjustment}
\label{sec:adjustment}

In the second scenario, the decision-maker is subject to a cost adjustment if the number of shutdowns differs from $N$. As in Scenario~1, there is a penalty if the number of shutoffs is above $N$. Conversely, if this number is below $N$, part of the anticipated value of load lost is recovered due to reduced outages, e.g., due to on-call personnel or avoided need for the extra fuel. We model the adjustment as proportional to the value of lost load for a day, $\lambda > 0$, and express it as: $\lambda \left( \sum_{t=1}^T u_{t} - N \right)$.

In Scenario 2, the total cost with adjustment is defined as:
\begin{align}
c_\text{a}\left(N, \uu \right) &= \sum_{t=1}^T \left( a_{t+1} u_{t} +  A_{t+1} f_{t+1}\left(\XX_{t+1}\right) (1-u_{t}) + \max\left\{{s_{1,t+1} \left( u_{t} - u_{t-1} \right), s_{2,t+1}\left(u_{t-1} - u_{t} \right) } \right\}\right)\nonumber\\
&\quad + \lambda \left(\sum_{t=1}^T u_{t} - N \right) + \overline{C}_I(N),
\label{eq:total_cost_sc2}
\end{align}
where the cost penalty in~\eqref{eq:total_cost} has been replaced by the cost adjustment.
The PSPS scheduling problem with cost adjustment is given by:
\begin{equation}
\begin{aligned}
& \min_{\substack{u_t\\ t=1,2,\ldots,T}} & & \E \left[ c_\text{a}\left(N, \uu\right) \right]\\
& \text{subject to}
& & u_t \in \left\{0, 1 \right\}\\
& & & u_{0} = u_{T+1} = 0.
\end{aligned}
\label{eq:min_unconstrained_1}
\end{equation}

The solution to Problem~\eqref{eq:min_unconstrained_1} is given by the following proposition. We later establish Proposition~\ref{prop:policy_psps_cost_adj}'s optimality in Section~\ref{ssec:analysis_adj}.
\begin{proposition}
Consider problem~\eqref{eq:min_unconstrained_1}.
On day $d = T - t$, after receiving the current observations, $\XX_t = \xx$, then $u_d =1 $ if
\begin{equation}
\E\left[\left.f_{d-1}\left(\XX_{d-1} \right)\right| \XX_{d} = \xx \right] \geq \frac{1}{A_{d-1}} \left(h_{d-1}(1,\xx) - h_{d-1}(0,\xx) + a_{d-1} + \lambda + \left(1-u_{d+1}\right) s_{{1},d-1} - u_{d+1} s_{{2},d-1}\right), \label{eq:policy_adj}
\end{equation}
where
\begin{align*}
h_{d}(u_{d+1} = u,\XX_{d+1} =\xx) &= \sum_{\xs \in \mathcal{X}} \max \left\{ h_{d-1}(1,\xs) +  \E\left[ \left.A_{d-1}f_{d-1}\left(\XX_{d-1}\right)\right| \XX_d = \xs \right] + u s_{{2},d-1},\right.\\
&\qquad\qquad\quad \left. h_{d-1}(0,\xs) + a_{d-1} + \lambda + (1-u)s_{{1},d-1} \right\}  \Pr\left[\left.\XX_d = \xs \right| \XX_{d+1} = \xx \right],
\end{align*}
with $h_{0}(u,\xx) = \sum_{\xs \in \XX}\E\left[ \left.A_{T+1}f_{T+1}\left(\XX_{T+1}\right)\right| \XX_{T} = \xs \right]\Pr\left[\left.\XX_{T} = \xs \right| \XX_{T-1} = \xx \right] + u s_{{2},T+1}$.
\label{prop:policy_psps_cost_adj}
\end{proposition}

Similarly to the first scenario, Proposition~\ref{prop:policy_psps_cost_adj} provides a closed-form optimal policy for the PSPS scheduling with cost adjustment problem that can be readily implemented by system operators. 

\subsection{Scenario~3: Minimum Number of PSPSs}
\label{sec:minPSPS}
In the third scenario, no budget for PSPSs has been set. Instead, we suppose that the system operator wishes to minimize the expected number of PSPS events such that the expected cost of the whole wildfire prevention strategy is below a threshold, which is itself less or equal to $C_I$. Similarly to previous sections, the strategy models both the initial, reduced investment and the costs of PSPSs. We suppose that a reduced investment $\tilde{C}_I < C_I$ was made to upgrade the power infrastructure and passively reduce wildfire risks. We let $\overline{\alpha} > 0$ be the difference between (i) the nominal investment $C_I$ and the (ii) the initial infrastructure investment $\tilde{C}_I$. The system operator must then schedule PSPSs such that the additional costs incurred are less or equal to $\overline{\alpha} = C_I - \tilde{C}_I$ and powerline-caused wildfire risks are mitigated in a cost-effective manner. We include Scenario~3 for completeness although it does not admit a readily implementable, closed-form policy.
We consider Scenario~3 because it may be of interest when the number of PSPSs is not constrained and difficult to forecast. Although we have not identified a readily implementable, closed-form policy at this time, we provide a value function-based analysis to optimally schedule PSPSs when its number is minimized.

Let $c_{t+1}:\left\{0,1\right\} \mapsto \RR$ be the operating cost function, i.e., the cost function without adjustment or penalty terms, for a single day $t+1$:
\begin{align*}
c_{t+1}\left(u_{t} \right) &= a_{t+1} u_{t} +  A_{t+1} f_{t+1}\left(\XX_{t+1}\right) (1-u_{t}) + \max\left\{{s_{1,t+1} \left( u_{t} - u_{t-1} \right), s_{2,t+1}\left(u_{t-1} - u_{t} \right) } \right\}.
\end{align*}
The minimum number of PSPS scheduling problem is:
\begin{equation}
\begin{aligned}
& \min_{\substack{u_t\\ t=1,2,\ldots,T}} & & \E \left[ \sum_{t=1}^T  u_t \right]\\
& \text{subject to}
& & u_t \in \left\{0, 1 \right\}\\
& & & u_{0} = u_{T+1} = 0 \\
& & & \E \left[ \sum_{t=1}^T  c_{t+1}\left(u_t \right) \right] \leq \overline{\alpha},
\label{eq:minPSPS}
\end{aligned}
\end{equation}
The optimal scheduling policy for problem~\eqref{eq:minPSPS} is presented in Proposition~\ref{prop:min_psps}.

\begin{proposition}
Consider the minimum number of PSPS scheduling problem~\eqref{eq:minPSPS} and assume it is feasible for $\XX_1 = \xx$. Let
\begin{align}
\alpha_{t} &= \begin{cases}
\overline{\alpha}, &\text{ if } t = 1\\
\phi^\star_{t}(\XX_t=\xx), &\text{ if }t \geq 2,
\end{cases}
\label{eq:policy_alpha}
\end{align}
where $\phi^\star_{t}(\XX_t=\xx)$ is defined by the following recursion
\begin{align}
\left(u^\star_t, \phi^\star_{t+1}\right) &= \begin{cases}
\left(1, \argmin_{\substack{\phi'(\xx') \in \Phi_{T-t}(\xx')}} \sum_{\xx' \in \mathcal{X}} \Pr[\left. \XX_{t+1}'=\xx' \right| \XX_{t}=\xx] V_{T-t-1}\left(\xx', \phi'(\xx') \right)\right),&\\
&\hspace{-9.80cm} \text{ if } \E\left[\left. c_{t+1}(0)\right| \XX_{t} = \xx\right] +\sum_{\xx' \in \mathcal{X}} {b}_{T-t-1}(\xx') \Pr[\left. \XX_{t+1}'=\xx' \right| \XX_{t}=\xx] > \alpha_t \\
&\hspace{-10.2cm}\text{ and } \E\left[\left. c_{t+1}(1)\right| \XX_{t} = \xx\right] +\sum_{\xx' \in \mathcal{X}} {b}_{T-t-1}(\xx') \Pr[\left. \XX_{t+1}'=\xx' \right| \XX_{t}=\xx] \leq \alpha_t \\
\left(0 , \argmin_{\phi'(\xx') \in \Phi_{T-t}(\xx')} \sum_{\xx' \in \mathcal{X}} \Pr[\left. \XX_{t+1}'=\xx' \right| \XX_{t}=\xx] V_{T-t-1}\left(\xx', \phi'(\xx') \right)\right),&\hspace{-0.25cm} \text{otherwise,} 
\end{cases}\label{eq:policy_u}
\end{align}
for $t=1,2,\ldots,T$ and where $\Phi_{T-t}(\xx,\alpha)$, $b_{T-t}(\xx,\alpha)$, and $V_T$ are respectively define in~\eqref{eq:b_tau},~\eqref{eq:phi_tau}, and~\eqref{eq:value_function_problem} from Section~\ref{ssec:value_func}.

Then, at day $t$, after receiving the current day's observations, $\XX_t = \xx$, it is optimal to declare a PSPS for the next day if $u_t^\star = 1$.
\label{prop:min_psps}
\end{proposition}

The derivation of Proposition~\ref{prop:min_psps}'s policy is given in Section~\ref{sec:analysis_scenario3}. We remark that while the problem formulation is of interest, to the author's best knowledge it does not admit a closed-form policy. This is a topic for future work.

\section{Analysis of Policies}
\label{sec:analysis}
In this section, we provide a comprehensive derivation for each optimal or asymptotically optimal policy. 

\subsection{Asymptotically Optimal Policy for Scenario~1}
\label{sec:analysis_costpenalty}
We now show that~\eqref{eq:min_unconstrained} is a relaxation of a modified multiple secretary problem~\citep{arlotto2019uniformly}, which has a closed-form optimal policy. Then, we establish the asymptotic exactness of the relaxation under certain conditions in Section~\ref{sec:relax}.

\subsubsection{Problem Reformulation.}
\label{sec:form}

In this section, we state the problem to which~\eqref{eq:min_unconstrained} is a relaxation.

We observe that for any $\gamma > 0$,~\eqref{eq:min_unconstrained} can be shown via the Lagrangian to be a relaxation of the following constrained problem: 
\begin{equation}
\begin{aligned}
& \min_{\substack{u_t\\ t=1,2,\ldots,T}} & & \E \left[\sum_{t=1}^T  a_{t+1} u_{t} +  A_{t+1} f_{t+1}\left(\XX_{t+1}\right) (1-u_{t}) + \max\left\{{s_{1,t+1} \left( u_{t} - u_{t-1} \right), s_{2,t+1}\left(u_{t-1} - u_{t} \right) } \right\}\right]\\
& \text{subject to}
& & u_t \in \left\{0, 1 \right\}\\
& & & u_{0} = u_{T+1} = 0 \\
& & & \E\left[\sum_{t=1}^T u_t \right] \leq N.
\end{aligned}
\label{eq:min_mag}
\end{equation}
In~\eqref{eq:min_mag}, the cost penalty for additional PSPS is now modelled as an expected budget constraint. 
We also consider~\eqref{eq:min_mag_noE} in which the expected constraint has been replaced by a deterministic constraint:
\begin{equation}
\begin{aligned}
& \min_{\substack{u_t\\ t=1,2,\ldots,T}} & & \E \left[\sum_{t=1}^T  a_{t+1} u_{t} +  A_{t+1} f_{t+1}\left(\XX_{t+1}\right) (1-u_{t}) + \max\left\{{s_{1,t+1} \left( u_{t} - u_{t-1} \right), s_{2,t+1}\left(u_{t-1} - u_{t} \right) } \right\}\right]\\
& \text{subject to}
& & u_t \in \left\{0, 1 \right\}\\
& & & u_{0} = u_{T+1} = 0 \\
& & & \sum_{t=1}^T u_t \leq N.
\end{aligned}
\label{eq:min_mag_noE}
\end{equation}
We remark that~\eqref{eq:min_mag} is a relaxation of~\eqref{eq:min_mag_noE} because the deterministic constraint implies that the expected one is satisfied. This follows from the fact that the~\eqref{eq:min_mag} entails that the budget constraint holds for all random variable sequence $\left\{\mathbf{X}_t \right\}_{t=1}^{T}$ outcomes. The differs from~\eqref{eq:min_mag_noE}, which can allow a PSPS budget above or below $N$ for some sequences as long as the budget constraint holds in the expected sense over all random variable sequences. We use the law of iterated expectations~\cite[Proposition 5.1]{ross2014first} to rewrite the objective of~\eqref{eq:min_mag_noE} as:
\begin{equation}
\begin{aligned}
& \min_{\substack{u_t\\ t=1,2,\ldots,T}} & & \E \Bigg[ \sum_{t=1}^T \E \big[  a_{t+1} u_{t} +  A_{t+1} f_{t+1}\left(\XX_{t+1}\right) (1-u_{t}) \\
& & & \qquad\qquad\quad\left.+ \max\left\{{s_{1,t+1} \left( u_{t} - u_{t-1} \right), s_{2,t+1}\left(u_{t-1} - u_{t} \right) } \right\} \right| \XX_{t} \big] \Bigg]\\
& \text{s.t.}
& & u_t \in \left\{0, 1 \right\}\\
& & & u_{0} = u_{T+1} = 0 \\
& & & \sum_{t=1}^T u_t \leq N.
\label{eq:min_mag_iterated}
\end{aligned}
\end{equation}
We re-write~\eqref{eq:min_mag_iterated} as an equivalent dynamic program as in~\citep{arlotto2019uniformly}. Let $w \in \RR$ be the cumulative cost incurred as of round $t$, i.e.,
\[
w = \sum_{i=1}^{t-1} a_{i+1} u_{i} +  A_{i+1} f_{i+1}\left(\xx_{i+1}\right) (1-u_{i}) + \max\left\{s_{{2},i+1}\left(u_{i-1} - u_{i} \right), s_{{1},i+1} \left( u_{i} - u_{i-1} \right)\right\},
\]
where $\xx_i$ is the realization of $\XX_i$ at round $i$. Let $d=T-t$ be the number of remaining days before the end of the PSPS program. To facilitate the analysis of the problem in its dynamic programing form, for the remainder of this subsection, all subscripts refer to the remaining number of decisions~$d$, e.g., $d=0$ or $d=-1$ is equivalent to $t=T$ and $t=T+1$, respectively.
Then, given decision $u_{d+1}$ and observations $\xx_{d+1}$, i.e., the decision and observations from the previous round, the value function is
\begin{equation}
\begin{aligned}
v_d\left(w,k\left|u,\xx\right.\right) = \sum_{\xs \in \mathcal{X}} \min \{ &v_{d-1}\left(\left.w + \E\left[\left. A_{d-1}f\left( \XX_{d-1} \right) \right| \XX_d = \xs \right] + u s_{{2},d-1},k\right|0,\xs\right),\\
& v_{d-1}\left(\left.w + a_{d-1} + (1-u)s_{{1}, d-1}, k-1 \right| 1 , \xs \right) \}\Pr\left[\left.\XX_d = \xs \right| \XX_{d+1} = \xx \right],
\end{aligned}
\label{eq:dp_psps_budget}
\end{equation}
with the boundary conditions at $d=0$ and $d-1=-1$:
\begin{align}
v_0(w,k|u,\xx) &= w + \E\left[\left. A_{-1}f\left( \XX_{-1} \right) \right| \XX_0 = \xs_j \right] + u s_{{2},-1},\label{eq:v_boundary}
\end{align}
$\text{ for all } w \in \RR, k=1, 2, \ldots, N, \xx \in \XX, u \in \left\{0,1 \right\}, \text{ and } d = 1, 2, \ldots, T$, and
\begin{align}
v_d(w,0|u,\xx) &= w + \sum_{i=0}^d\sum_{\xs \in \mathcal{X}}\E\left[ \left.A_{i-1}f_{i-1}\left(\XX_{i-1}\right)\right| \XX_{i} = \xs \right]\mathbf{P}^{(i+1)}_{\xx,\xs} + u s_{{2},d-1},
\label{eq:v_0}
\end{align}
for all $w \in \RR, \xx \in \XX, u \in \left\{0,1 \right\}, \text{ and } d = 1, 2, \ldots, T$.

The dynamic program~\eqref{eq:dp_psps_budget}$-$\eqref{eq:v_0} admits a closed-form, easily-implementable optimal policy. We then show that problem~\eqref{eq:min_mag_noE}, which is equivalent to~\eqref{eq:dp_psps_budget}$-$\eqref{eq:v_0}, is asymptotically equivalent to the original problem~\eqref{eq:min_unconstrained} under certain conditions. Problem~\eqref{eq:dp_psps_budget}$-$\eqref{eq:v_0} can be interpreted as an instance of the multiple secretary problem~\citep{arlotto2019uniformly} with the addition of a switching cost and a regularizer. The switching cost further adds a time-dependent constraint to the problem. This will lead to more conservative decision-making in general but to a more sensitive policy immediately following a PSPS. The regularizer takes the form of an offset which penalizes the worst decision that can be taken, i.e., selecting a day with low expected wildfire risks. It can be interpreted as a sparsity regularizer. Consequently, it may lead to decision sequences for which the budget may not be fully used, which was not considered in prior multiple secretary formulations. We adapt the results~\citep[Appendix C]{arlotto2019uniformly} to include switching costs and the regularizer and obtain the optimal policy for~\eqref{eq:dp_psps_budget}$-$\eqref{eq:v_0}, as shown above in Proposition~\ref{prop:policy_psps}. 

\proof{Proof of Proposition~\ref{prop:policy_psps}.}
Let $e_{d-1}(\xs) = \E\left[ \left.A_{d-1}f_{d-1}\left(\XX_{d-1}\right)\right| \XX_{d} = \xs \right]$ and $p_\xx(\xs) = \Pr\left[\left.\XX_d = \xs \right| \XX_{d+1} = \xx \right]$ to simplify notation.
We first show that
\begin{equation}
v_{d}\left(w,k\left|u,\xx\right.\right) = w + g_{d}(k|u,\xx),
\label{eq:recursion}
\end{equation}
for $d=0, 1, \ldots, T$. We proceed by induction. In the base case $d=0$,~\eqref{eq:recursion} holds trivially because of the boundary conditions and for $d =1$, we have two instances: $k=0$ and $k\geq 1$. For $k=0$,
\begin{align}
v_1(w,0|u,\xx) &= w + \sum_{\xs \in \XX}\E\left[ \left.A_{0}f_{0}\left(\XX_{0}\right)\right| \XX_{1} = \xs \right]\Pr\left[\left.\XX_1 = \xs \right| \XX_{2} = \xx \right]+ u s_{{2},0}. \label{v1_k0}
\end{align}
For $k\geq 1$,
\begin{align}
v_1(w,k|u,\xx) =& \sum_{\xs \in \mathcal{X}} \min \{ v_0\left(\left.w + e_{0}(\xs) + u s_{{2},0},k\right|0,\xs\right), v_0\left(\left.w + a_{0} + (1-u)s_{{1}, 0}, k-1 \right| 1 , \xs \right) \} p_\xx(\xs) \nonumber\\
=& \sum_{\xs \in \mathcal{X}} \min \{ w + e_{0}(\xs) + u s_{{2},0}, w + a_{0} + (1-u)s_{{1}, 0} + s_{{2},-1} \} p_\xx(\xs) \nonumber\\
=& w + \sum_{\xs \in \mathcal{X}} \min \{ e_{0}(\xs) + u s_{{2},0}, a_{0} + (1-u)s_{{1}, 0} + s_{{1},-1} \} p_\xx(\xs) \label{eq:v1_k}
\end{align}
where the second equality follows from the boundary condition~\eqref{eq:v_boundary}. Similarly, for $g_1(k|u,\xx)$, we have for $k=0$:
\begin{align}
g_{1}(0|u,\xx) &=\sum_{\xs \in \XX}\E\left[ \left.A_{0}f_{0}\left(\XX_{0}\right)\right| \XX_{1} = \xs \right]\Pr\left[\left.\XX_1 = \xs \right| \XX_{2} = \xx \right] + u s_{{2},0} \label{eq:g1_k0}.
\end{align}
For $k \geq 1$, we have
\begin{align}
g_{1}(k|u,\xx) &= \sum_{\xs \in \mathcal{X}} \max \left\{ g_{0}(k|0,\xs) +  e_0(\xs) + u s_{{2},0}, g_{0}(k-1|1,\xs) + a_{0} + (1-u) s_{{1},0}\right\}  p_\xx(\xs)\nonumber\\
&= \sum_{\xs \in \mathcal{X}} \max \left\{ e_0(\xs) + u s_{{2},0}, a_{0} + (1-u) s_{{1},0} + s_{{2},-1}\right\}  p_\xx(\xs), \label{eq:g1_k}
\end{align}
where we last used the boundary conditions~\eqref{g_0d}. Using~\eqref{eq:g1_k0} in~\eqref{v1_k0} and~\eqref{eq:v1_k} in~\eqref{eq:g1_k}, we obtain
\[
v_1\left(w,k\left|u,\xx\right.\right) = w + g_{1}(k|u,\xx),
\]
which establishes the base case. Now, let $d \to d +1$. Then for $k=0$, by definition we have
\begin{align}
v_{d+1}(w,0|u,\xx) &= w + \sum_{i=0}^{d+1}\sum_{\xs \in \mathcal{X}}\E\left[ \left.A_{i-1}f_{i-1}\left(\XX_{i-1}\right)\right| \XX_{i} = \xs \right]\mathbf{P}^{(i+1)}_{\xx,\xs} + u s_{{2},d} \label{eq:vdk0}\\
g_{d+1}(0|u,\xx) &=\sum_{i=0}^d\sum_{\xs \in \mathcal{X}}\E\left[ \left.A_{i-1}f_{i-1}\left(\XX_{i-1}\right)\right| \XX_{i} = \xs \right]\mathbf{P}^{(i+1)}_{\xx,\xs} + u s_{{2},d}. \label{eq:gdk0}
\end{align}
For $k\geq1$, we obtain
\begin{align*}
v_{d+1}\left(w,k\left|u,\xx\right.\right) &= \sum_{\xs \in \mathcal{X}} \min \{ v_{d}\left(\left.w + e_{d}(\xs) + u s_{{2},d},k\right|0,\xs\right), v_{d}\left(\left.w + a_{d} + (1-u)s_{{1}, d}, k-1 \right| 1 , \xs \right) \}p_\xx(\xs)\\
&= w+ \sum_{\xs \in \mathcal{X}} \min \{ v_{d}\left(\left.w + e_{d}(j) + u s_{{2},d},k\right|0,\xs\right) + e_{d}(\xs) + u s_{{2},d} -w  - e_{d}(\xs) - u s_{{2},d},\\
&\quad v_{d}\left(\left.w + a_{d} + (1-u)s_{{1}, d}, k-1 \right| 1 , \xs \right) + a_{d} + (1-u)s_{{1}, d} - w - a_{d} - (1-u)s_{{1}, d} \}p_\xx(\xs)\\
&= w+ \sum_{j \in \mathcal{X}} \min \{ g_{d}\left(\left.k\right|0,\xs\right) + e_{d}(\xs) + u s_{{2},d}, g_{d}\left(\left. k-1 \right| 1 , \xs \right) + a_{d} + (1-u)s_{{1}, d} \}p_\xx(\xs)\\
&= w+ g_{d+1}\left(\left.k\right|u,\xx\right),
\end{align*}
where we used the induction hypothesis~\eqref{eq:recursion} and the definition of $g$ given in~\eqref{eq:def_g} to obtain the third and final equality, respectively. This completes the induction proof and we have established~\eqref{eq:recursion} for $d=0,1,\ldots, T$.

We now derive an optimal policy using~\eqref{eq:def_g}. The structure of $v_d\left(w,k\left|u_{d+1},\XX_{d+1}\right.\right)$ implies that $u_d = 1$ if and only if:
\begin{align*}
v_{d-1}\left(\left.w + \E\left[\left. A_{d-1}f\left( \XX_{d-1} \right) \right| \XX_d = \xx \right] + u_{d+1} s_{{2},d-1},k\right|0,\xx\right) \geq&\\
&\hspace{-0.5cm}v_{d-1}\left(\left.w + a_{d-1} + (1-u_{d+1})s_{{1}, d-1}, k-1 \right| 1 , \xx \right).
\end{align*}
Substituting~\eqref{eq:recursion} yields
\begin{align*}
g_{d-1}\left(\left.k\right|0,\xx\right) + w + \E\left[\left. A_{d-1}f\left( \XX_{d-1} \right) \right| \XX_d = \xx \right] + u_{d+1} s_{{2},d-1} \geq g_{d-1}\left(\left.k-1\right|1,\xx\right)&\\
&\hspace{-1cm}+ w + a_{d-1} + (1-u_{d+1})s_{{1}, d-1}.
\end{align*}
Thus, we have $u_d = 1$ if and only if
\[
\E\left[\left. f\left( \XX_{d-1} \right) \right| \XX_d = \xx \right] \geq \frac{1}{A_{d-1}}\left( g_{d-1}\left(\left.k-1\right|1,\xx\right) - g_{d-1}\left(\left.k\right|0,\xx\right) + a_{d-1} + (1-u_{d+1})s_{{1},d-1} - u_{d+1} s_{{2},d-1}\right),
\]
which completes the proof. \hfill \Halmos
\endproof

In the above policy, the difference $g_{d-1}\left(\left.k\right|0,\xx\right) - g_{d-1}\left(\left.k-1\right|1,\xx\right)$ ensures a PSPS is scheduled only if the cost of the wildfire at the next round weighted by the expected risk is greater than the expected cost reduction of a PSPS when $k$ have been previously called.
This is done to account for the budget constraint. The other terms modify the threshold to account for fixed costs of calling PSPS. In particular, the switching costs promotes consecutive over isolated de-energization events. The regularizer integrates the fact that the optimal strategy may be to call less than $N$ PSPSs. This differs from~\citep{arlotto2019uniformly} in which the optimal solution always uses the full budget.

\subsubsection{Asymptotically Exact Relaxation.}
\label{sec:relax}

We now show that under certain conditions, the average cost or per-round form of~\eqref{eq:min_mag_noE}, in which the total expected cost is scaled by the time horizon, is asymptotically equivalent to the original problem~\eqref{eq:min_unconstrained}'s per-round form, i.e., the latter is an asymptotically exact relaxation of the former. 

We first establish sufficient conditions for~\eqref{eq:min_mag} to be an exact relaxation of~\eqref{eq:min_unconstrained} in Lemma~\ref{lem:lag_exact}. To this end, let $r:\left\{0,1 \right\}^T \mapsto \RR$ be the objective function of~\eqref{eq:min_mag_noE}. Let $P:\left\{0,1 \right\}^T \mapsto \RR^+$ where $P(\uu) = \E\left[\max\left\{0,\sum_{t=1}^T u_t - N \right\}\right]$ be the penalty for additional PSPSs. Recall that $\gamma$ is the coefficient of the penalty term in~\eqref{eq:min_unconstrained}.
\begin{lemma}
Suppose $N < T$. Let
\begin{align*}
\alpha =  \min\left\{  r(\uu) \left| \uu \in \left\{0 ,1 \right\}^T, u_0=u_{T+1}=0, \E\left[\sum_{t=1}^T u_t \right] \leq N \right.\right\}&\\
&\hspace{-0.5cm}- \min\left\{ r(\uu) \left| \uu \in \left\{0 ,1 \right\}^T, u_0=u_{T+1}=0 \right.\right\}.
\end{align*}
Then for all $\gamma > \alpha$,~\eqref{eq:min_unconstrained} is an exact relaxation of~\eqref{eq:min_mag}.
\label{lem:lag_exact}
\end{lemma}

\proof{Proof of Lemma~\ref{lem:lag_exact}.} 
We first observe that: 
$$\min\left\{ P(\uu) \left| \uu \in \left\{0 ,1 \right\}^T, u_0=u_{T+1}=0, \E\left[\sum_{t=1}^T u_t \right] > N \right.\right\} = 1.$$
Then by~\citep[Theorem 1]{sinclair1986exact}, for $\gamma > \alpha$, the optima of~\eqref{eq:min_unconstrained} are also optima of~\eqref{eq:min_mag}.

Second, we show the converse which was omitted in~\citep{sinclair1986exact}, i.e., the optima of~\eqref{eq:min_mag} are optimal for~\eqref{eq:min_unconstrained} when $\gamma > \alpha$. Let $o_\eqref{eq:min_unconstrained}$ and $o_\eqref{eq:min_mag}$ be the optimal value of~\eqref{eq:min_unconstrained} and~\eqref{eq:min_mag}, respectively. Let $\uu^\star$ be an optimum of~\eqref{eq:min_mag}. For all $\gamma > 0$, we have
\begin{align*}
o_\eqref{eq:min_mag} &= r(\uu^\star)\\
&= r(\uu^\star) + \gamma P(\uu^\star),
\end{align*}
because $\uu^\star$ is feasible for~\eqref{eq:min_mag}, i.e., satisfies $\E\left[\sum_{t=1}^T u_t \right] \leq N$ and, thus, $P(\uu^\star) = 0$. By~\citep[Theorem 1]{sinclair1986exact}, $o_\eqref{eq:min_unconstrained} = o_\eqref{eq:min_mag}$ for $\gamma > \alpha$. Hence,
\[
r(\uu^\star) + \gamma P(\uu^\star) = o_\eqref{eq:min_unconstrained},
\]
and $\uu^\star$ is also an optimum of~\eqref{eq:min_unconstrained}. Therefore, if $\gamma > \alpha$, then~\eqref{eq:min_unconstrained} and~\eqref{eq:min_mag} have the same optima. It follows that~\eqref{eq:min_unconstrained} is an exact relaxation of~\eqref{eq:min_mag} under this condition.\hfill \Halmos
\endproof

Lemma~\ref{lem:lag_exact} states that if the marginal cost of being over the PSPS budget is large enough, then the relaxation is exact.

We now show that~\eqref{eq:min_mag_noE}'s and~\eqref{eq:min_mag}'s per-round forms are equivalent, almost surely, when the number of days and the budget constraint threshold $N$ grow to infinity with a fixed ratio $\rho \in (0,1)$. 

Consider the Markov decision processes $\left(\mathcal{X}, \left\{0,1\right\}, \mathbf{P}, r_t \right)$ used in~\eqref{eq:min_mag} and~\eqref{eq:min_mag_noE}. Then, the Markov chain induced by any policy is the same for both MDPs because $\XX_t$, $t=1,2,\ldots,T$, models natural phenomena and its evolution is independent of the decisions $\mathbf{u}$. Let $\mathcal{M}_\XX$ denote this Markov chain.

Recall that $\mathcal{M}_\XX$ is ergodic if it is recurrent and aperiodic. An ergodic Markov chain $\mathcal{M}_\XX$ possesses a unique stationary distribution $\mathbf{s} \in \RR^{n}_+$ such that $\mathbf{s} = \mathbf{s} \mathbf{P}$ and $\mathbf{1}^\top \mathbf{s} = 1$. Finally, $\mathcal{M}_\XX$ is stationary if the initial state of the Markov chain is distributed according to $\mathbf{s}$, i.e., $\XX_0 \sim \mathbf{s}$. For $\mathcal{M}_\XX$ stationary, we therefore have $\XX_t \overset{\text{d}}{=} \XX \overset{\text{d}}{=} \XX_0 \sim \mathbf{s}$ for all $t \geq 0$, where $\overset{\text{d}}{=}$ means that they are equal in distribution.

Let $\Pi$ be the policy space. Let $\Pi^{\eqref{eq:min_mag_noE}}$ be the set of all feasible policies for Problem~\eqref{eq:min_mag_noE}. We now present a lemma stating an equivalent average per-round constraint for policies in $\Pi^{\eqref{eq:min_mag_noE}}$.

\begin{lemma}
Suppose $\mathcal{M}_X$ is ergodic and stationary, and consider the policy $\pi^{\eqref{eq:min_mag_noE}} \in \Pi$. 
Let $T, N \to \infty$, where $T$ is an integer and $N=\rho T$. Then $\pi^{\eqref{eq:min_mag_noE}} \in \Pi^{\eqref{eq:min_mag_noE}}$ if and only if $\E\left[ \pi^{\eqref{eq:min_mag_noE}}(\XX)\right] \leq \rho$ almost surely.
\label{lem:per_round}
\end{lemma}

\proof{Proof of Lemma~\ref{lem:per_round}.}
Let $\pi^{\eqref{eq:min_mag_noE}} \in \Pi^{\eqref{eq:min_mag_noE}}$. By definition, we have for all $T, N \geq 1$,
\[
\sum_{t=1}^T u_t^{\pi^{\eqref{eq:min_mag_noE}}} \leq N,
\]
and therefore
\[
\frac{1}{T}\sum_{t=1}^T \pi^{\eqref{eq:min_mag_noE}}(\XX_t) \leq \frac{N}{T},
\]
where we now consider the per-round average constraint and used the fact that $u_t^{\pi^{\eqref{eq:min_mag_noE}}} = \pi^{\eqref{eq:min_mag_noE}}(\XX_t)$.

We take the limit as $T,N \to \infty$, with $T \in \mathbb{N}$ and $N/T$ kept constant at $\rho$. We have
\begin{align}
\lim_{\substack{T,N \to +\infty \\ \text{s.t.} T \in \mathbb{N} \\ N = \rho T}}\frac{1}{T}\sum_{t=1}^T \pi^{\eqref{eq:min_mag_noE}}(\XX_t) &\leq \lim_{\substack{T,N \to +\infty \\ \text{s.t.} T \in \mathbb{N} \\ N = \rho T}} \frac{N}{T},\nonumber\\
&= \rho \label{eq:per_round_ratio}
\end{align}
Using a strong law of large number for ergodic Markov chains~\citep[Chapter 1, Theorem 74]{serfozo2009basics} we obtain for the left-hand side:
\begin{align*}
\lim_{\substack{T,N \to +\infty \\ \text{s.t.} N = \rho T}}\frac{1}{T}\sum_{t=1}^T \pi^{\eqref{eq:min_mag_noE}}(\XX_t) = \sum_{\xx \in \mathcal{X}} \pi^{\eqref{eq:min_mag_noE}}(\xx) \mathbf{s}(\xx) \quad \text{ almost surely},
\end{align*}
if the right-hand term is absolutely convergent. This condition is met in~\eqref{eq:min_mag_noE} because $\pi^{\eqref{eq:min_mag_noE}}(\XX) \in \left\{0,1\right\}$ for all $\xx \in \mathcal{X}$ and $\card \mathcal{X}$ is finite. We rewrite~\eqref{eq:per_round_ratio} as
\begin{equation}
\E_{\XX \sim \mathbf{s}}\left[ \pi^{\eqref{eq:min_mag_noE}}(\XX) \right] \leq \rho \quad \text{ almost surely}. \label{eq:per_round_exp}
\end{equation}
Hence, if $T,N \to \infty$ such that $T \in \mathbb{N}$ and  $\rho = N/T$, then the constraint defining $\Pi^{\eqref{eq:min_mag_noE}}$ is equivalent to~\eqref{eq:per_round_exp}. It follows that $\pi^{\eqref{eq:min_mag_noE}} \in \Pi^{\eqref{eq:min_mag_noE}}$ if and only if $\E\left[ \pi^{\eqref{eq:min_mag_noE}}(\XX)\right] \leq \rho$. \hfill \Halmos
\endproof

Similarly, we now present the per-round constraint that defines all policies in $\Pi^{\eqref{eq:min_mag}}$.

\begin{lemma}
Suppose $\mathcal{M}_X$ is ergodic and stationary and let $\pi^{\eqref{eq:min_mag}} \in \Pi$. As $T, N\to + \infty$, where $T$ is an integer and $N=\rho T$, $\pi^{\eqref{eq:min_mag}} \in \Pi^{\eqref{eq:min_mag}}$ if and only if $\E\left[\pi^{\eqref{eq:min_mag}}(\XX) \right] \leq \rho$.
\label{lem:per_round_bis}
\end{lemma}

\proof{Proof of Lemma~\ref{lem:per_round_bis}.}
The first steps are similar to Lemma~\ref{lem:per_round}'s proof. Let $\Pi^{\eqref{eq:min_mag}}$ be the set of all feasible policies for~\eqref{eq:min_mag}. Let $\pi^{\eqref{eq:min_mag}} \in \Pi^{\eqref{eq:min_mag}}$. Be definition, 
\begin{align}
\E\left[\sum_{t=1}^T u_t^{\pi^{\eqref{eq:min_mag}}} \right] &\leq N \label{eq:ori_constraint}\\
\iff \sum_{t=1}^T \E\left[u_t^{\pi^{\eqref{eq:min_mag}}} \right] &\leq N \nonumber\\
\iff \frac{1}{T}\sum_{t=1}^T \E\left[u_t^{\pi^{\eqref{eq:min_mag}}} \right] &\leq \frac{N}{T}\nonumber\\
\iff \frac{1}{T}\sum_{t=1}^T \E\left[\pi^{\eqref{eq:min_mag}}(\XX_t) \right] &\leq \frac{N}{T}\nonumber
\end{align}
We evaluate the limit on both sides as $T,N \to +\infty$ when $T \in \mathbb{N}$ and the ratio $\rho=N/T$ is kept constant. We have
\begin{align}
\lim_{\substack{T,N \to +\infty \\ \text{s.t.} N = \rho T}} \frac{1}{T}\sum_{t=1}^T \E\left[\pi^{\eqref{eq:min_mag}}(\XX_t) \right] &\leq \lim_{\substack{T,N \to +\infty \\ \text{s.t.} N = \rho T}} \frac{N}{T}\nonumber\\
& = \rho. \label{eq:stat_rho}
\end{align}
Because $\mathcal{M}_\XX$ is stationary, we have $\XX_t \overset{d}{=} \XX \sim \mathbf{s}$ for all $t \geq 0$ and thus
\begin{align}
\lim_{\substack{T,N \to +\infty \\ \text{s.t.} N = \rho T}} \frac{1}{T}\sum_{t=1}^T \E\left[\pi^{\eqref{eq:min_mag}}(\XX_t) \right] &= \lim_{\substack{T,N \to +\infty \\ \text{s.t.} N = \rho T}} \frac{1}{T}\sum_{t=1}^T \E\left[\pi^{\eqref{eq:min_mag}}(\XX) \right]\nonumber\\
& = \E\left[\pi^{\eqref{eq:min_mag}}(\XX) \right]. \label{eq:stat_exp}
\end{align}
Using~\eqref{eq:stat_rho} and~\eqref{eq:stat_exp} leads to
\begin{equation}
\E\left[\pi^{\eqref{eq:min_mag}}(\XX) \right] \leq \rho,
\label{eq:stat_condition}
\end{equation}
and we have showed that~\eqref{eq:ori_constraint} is equivalent to~\eqref{eq:stat_condition}. It follows that $\pi^{\eqref{eq:min_mag}} \in \Pi^{\eqref{eq:min_mag}}$ if and only if $\E\left[\pi^{\eqref{eq:min_mag}}(\XX) \right] \leq \rho$ as $T \to \infty$ and $N = \rho T$.\hfill \Halmos
\endproof

Next, we establish the asymptotic relationship between round-averaged forms of~\eqref{eq:min_mag} and~\eqref{eq:min_mag_noE}.

\begin{theorem}
Suppose $\mathcal{M}_\XX$ is ergodic and stationary. The per-round average form of~\eqref{eq:min_mag} is an asymptotically exact relaxation of~\eqref{eq:min_mag_noE}'s round-averaged formulation, almost surely, as $T,N \to \infty$ with $\rho = N/T$, $T \in \mathbb{N}$. 
\label{th:equivalent}
\end{theorem}

\proof{Proof of Theorem~\ref{th:equivalent}.}
From Lemmas~\ref{lem:per_round} and~\ref{lem:per_round_bis}, we have $\pi^{\eqref{eq:min_mag_noE}} \in \Pi^{\eqref{eq:min_mag_noE}}$ if and only if $\E\left[\pi^{\eqref{eq:min_mag_noE}}(\XX) \right] \leq \rho$ almost surely and $\pi^{\eqref{eq:min_mag}} \in \Pi^{\eqref{eq:min_mag}}$ if and only if $\E\left[\pi^{\eqref{eq:min_mag}}(\XX) \right] \leq \rho$ when $T,N \to \infty$ with $\rho = N/T$. Therefore, $\Pi^{\eqref{eq:min_mag}} = \Pi^{\eqref{eq:min_mag_noE}}$ almost surely as $T,N \to \infty$ with $\rho = N/T$. This implies that the per-round average form of problems~\eqref{eq:min_mag} and~\eqref{eq:min_mag_noE} are asymptotically equivalent almost surely because their objective functions are the same. In other words, the per-round form of~\eqref{eq:min_mag} is an asymptotically exact relaxation of the~\eqref{eq:min_mag_noE}'s per-round form.\hfill \Halmos
\endproof

Lastly, we relate the original problem~\eqref{eq:min_unconstrained} to the final form~\eqref{eq:min_mag_noE} for which we have an optimal PSPS scheduling policy.

\begin{theorem}
If $\gamma > \alpha$ and $\mathcal{M}_\XX$ is ergodic and stationary, then, the per-round form of~\eqref{eq:min_unconstrained} is an exact relaxation of~\eqref{eq:min_mag_noE}'s, almost surely, as $T, N \to \infty$ where $T$ is an integer and $N=\rho T$. Otherwise, the expected minimum cumulative loss of~\eqref{eq:min_mag_noE} is greater than or equal to the minimum of~\eqref{eq:min_unconstrained}.
\label{th:policy_and_original_prob}
\end{theorem}

\proof{Proof of Theorem~\ref{th:policy_and_original_prob}.}
We remark that the averaged problem, i.e., dividing the objective function by $T$, only scales down the objective, which remains finite as $T\to\infty$.  First, we prove the asymptotic exactness of the relaxation when $\gamma > \alpha$ and $\mathcal{M}_\XX$ is stationary and as $T, N \to \infty$ such that $N=\rho T$.
Problem~\eqref{eq:min_unconstrained} is a relaxation of~\eqref{eq:min_mag}. Because we have assumed $\gamma > \alpha$, Lemma~\ref{lem:lag_exact} ensures that the relaxation is exact. Thus,~\eqref{eq:min_mag}'s optima are also optimal for~\eqref{eq:min_unconstrained}. Problem~\eqref{eq:min_mag} is in turn a relaxation of~\eqref{eq:min_mag_noE}. Under the theorem's assumption, we can invoke Theorem~\ref{th:equivalent}, which establishes that the per-round form of~\eqref{eq:min_mag} is an asymptotically exact relaxation of the per-round problem~\eqref{eq:min_mag_noE}. Similarly,~\eqref{eq:min_mag_noE}'s optima are asymptotically optimal for~\eqref{eq:min_mag} almost surely as $T, N \to \infty$ such that $N=\rho T$. The optima of~\eqref{eq:min_mag_noE} are, therefore, asymptotically optimal for~\eqref{eq:min_unconstrained} almost surely as well when $\gamma > \alpha$. The relaxation is, therefore, exact asymptotically.

Finally, we discuss the case for which the assumptions are not satisfied. From the above justification, we have that regardless of the values of $\gamma$, $N$ and $T$, and the stationarity of $\mathcal{M}_\XX$,~\eqref{eq:min_unconstrained} is a relaxation of~\eqref{eq:min_mag_noE}. Problem~\eqref{eq:min_mag_noE}'s optimal cost are, therefore, either greater than or equal to the minimum of~\eqref{eq:min_unconstrained}.\hfill \Halmos
\endproof

We note that under any conditions, the decisions provided by Proposition~\ref{prop:policy_psps}'s policy are always feasible for PSPS scheduling and lead to a minimum that is equal to or greater than~\eqref{eq:min_unconstrained}. This is because~\eqref{eq:dp_psps_budget} is a reformulation of~\eqref{eq:min_unconstrained} for which Proposition~\ref{prop:policy_psps} provides the optimal solution.

\subsection{Optimal Policy for Scenario~2}
\label{ssec:analysis_adj}
We now prove Proposition~\ref{prop:policy_psps_cost_adj}.
We remove $\overline{C}_I(N)$ and $-\lambda N$ from the problem formulation because they are constants and do not impact the minima. We write~\eqref{eq:min_unconstrained_1} as
\begin{equation}
\begin{aligned}
& \min_{\substack{u_t\\ t=1,2,\ldots,T}} & & \E \Bigg[\sum_{t=1}^T  a_{t+1} u_{t} +  A_{t+1} f_{t+1}\left(\XX_{t+1}\right) (1-u_{t}) + \lambda u_t\\
& & & \qquad \qquad+ \max\left\{s_{{2},t+1}\left(u_{t-1} - u_{t} \right), s_{{1},t+1} \left( u_{t} - u_{t-1} \right) \right\}\Bigg]\\
& \text{subject to}
& & u_t \in \left\{0, 1 \right\}\\
& & & u_{0} = u_{T+1} = 0.
\end{aligned}
\label{eq:min_mag_1}
\end{equation}
Using the law of iterated expectation~\cite[Proposition 5.1]{ross2014first}, as in~\eqref{eq:min_mag_iterated}, gives
\begin{equation}
\begin{aligned}
& \min_{\substack{u_t\\ t=1,2,\ldots,T}} & & \E \Bigg[ \sum_{t=1}^T \E \big[  a_{t+1} u_{t} +  A_{t+1} f_{t+1}\left(\XX_{t+1}\right) (1-u_{t}) + \lambda u_t\\
& & & \qquad \qquad  \qquad+ \max\left\{s_{{2},t+1}\left(u_{t-1} - u_{t} \right), s_{{1},t+1} \left( u_{t} - u_{t-1} \right) \right\} \big| \XX_{t} \big] \Bigg]\\
& \text{subject to}
& & u_t \in \left\{0, 1 \right\}\\
& & & u_{0} = u_{T+1} = 0.
\end{aligned}
\label{eq:min_mag_iterated_1}
\end{equation}
Problem~\eqref{eq:min_mag_iterated_1} represents a Markov decision process, in which a decision maker observes the random variable $\XX_t$ at time $t$ and then make the decision, $u_t$, to call a PSPS for the next day. Similarly to Section~\ref{sec:analysis_costpenalty}, let $d=T-t$ and let all subscripts now refer to the backward time index $d$. Let $z_d\left(w\left|u_{d+1},\XX_{d+1}\right.\right)$ be the expected total cost given the past decision $u_{d+1}$ and the observation vector $\XX_{d+1}$, when $d$ days remains and the accumulated cost is $w$. Problem~\eqref{eq:min_mag_iterated_1} can be expressed as the following dynamic program:
\begin{align}
z_d\left(w\left|u,\xx\right.\right) = \sum_{\xs \in \mathcal{X}} \min \{ &z_{d-1}\left(\left.w + \E\left[\left. A_{d-1}f\left( \XX_{d-1} \right) \right| \XX_d = \xs \right] + u s_{{2},d-1}\right|0,\xs\right), \nonumber\\
& z_{d-1}\left(\left.w + a_{d-1} + \lambda + (1-u)s_{{1}, d-1}\right| 1 , \xs \right) \}\Pr\left[\left.\XX_d = \xs \right| \XX_{d+1} = \xx \right],
\label{eq:dp_adj}
\end{align}
with the boundary conditions: 
\[
z_0\left(w\left|u,\xx\right.\right) = w +\sum_{\xs \in \mathcal{X}}\E\left[ \left.A_{-1}f_{-1}\left(\XX_{-1}\right)\right| \XX_{0} = \xs \right]\Pr\left[\left.\XX_0 = \xs \right| \XX_{1} = \xx \right] + u s_{{2},-1},
\]
for all $w \in \RR, \xx \in \mathcal{X}$, and $u \in \{0,1 \}$.
The optimal policy for the cost adjustment problem is given in Proposition~\ref{prop:policy_psps_cost_adj} of Section~\ref{sec:adjustment}. Its proof is stated next.

\proof{Proof of Proposition~\ref{prop:policy_psps_cost_adj}.}
Proposition~\ref{prop:policy_psps_cost_adj} follows from the same proof technique as for Proposition~\ref{prop:policy_psps}, in which we use the relation:
\[
z_d(w\left| u, \xx\right. ) = w + h_d(u, \xx)
\]
for all $w \in \RR$, $u \in \left\{0,1\right\}$, $\xx \in \XX$ and $d=T, T-1, \ldots, 1, 0$. Then, given the structure of $z_d(w, u, \xx )$, the optimal policy is given by a threshold policy as well. Using the recursion and solving for $\E\left[\left.f_{d-1}\left(\XX_{d-1} \right)\right| \XX_{d} = \xx \right]$ yields~\eqref{eq:policy_adj}.\hfill \Halmos
\endproof

\subsection{Optimal Policy for Scenario~3}
\label{sec:analysis_scenario3}
We discuss Scenario~3 in which the number of PSPSs is minimized. Problem~\eqref{eq:minPSPS} can be re-written as
\begin{equation}
\begin{aligned}
& \min_{\substack{u_t\\ t=1,2,\ldots,T}} & & \E \left[ \sum_{t=1}^T \E \left[ \left.  u_t \right| \XX_{t} \right] \right]\\
& \text{subject to}
& & u_t \in \left\{0, 1 \right\}\\
& & & u_{0} = u_{T+1} = 0 \\
& & & \E \left[ \sum_{t=1}^T \E \left[ \left.  c_{t+1}\left(u_t \right) \right| \XX_{t} \right] \right] \leq \overline{\alpha},
\label{eq:minPSPS_2}
\end{aligned}
\end{equation}
where we have used the law of iterated expectations~\cite[Proposition 5.1]{ross2014first} to re-express the objective and the constraint functions. Problem~\eqref{eq:minPSPS_2} can be solved using value iteration as in~\citep{chen2004dynamic} and computing the optimal policy as in~\citep{chen2007non}. Adopting~\citep{chen2004dynamic}'s notation, we let $\Phi_\tau\left(\xx\right)$ be the set of feasible constraint thresholds for a problem starting at day $t$, defined as
\begin{equation}
\Phi_\tau(\XX_t = \xx) = \left\{\alpha \in \RR \left| \ \exists \right. \uu \in \left\{ 0, 1 \right\}^\tau, u_0=u_{
T+1}=0 \text{ such that }  \sum_{i=t}^{t+\tau} \E \left[ \left.  c_{i+1}\left(u_i \right) \right| \XX_{t} \right] \leq \alpha\right\}. \label{eq:PhiT}
\end{equation}
Let $\phi: \mathcal{X} \mapsto \RR^n$ be the constraint threshold for the initial state $\xx \in \mathcal{X}$. Let $\overline{V} \geq 0$ be a large scalar. Let $V_\tau(\xx,\alpha)$ be the value function or total expected cost for problem~\eqref{eq:minPSPS_2} with the time horizon $\tau$, initial state $\xx$, and constraint threshold $\alpha$. 
Based on~\citep[Theorem 3.1]{chen2004dynamic}, we re-express~\eqref{eq:minPSPS_2} as the dynamic program
\begin{equation}
V_{\tau+1}(\xx,\alpha) = \mathcal{T}_{F_\tau(\xx, \alpha)} V_\tau(\xx,\alpha),
\label{eq:dp_minPSPS}
\end{equation}
for $\tau \geq 1$, where the operator $\mathcal{T}_{F_\tau}$ is defined as:
\begin{equation}
\mathcal{T}_{F_\tau(\xx, \alpha)} V(\xx,\alpha) = \begin{cases}
\min_{(u, \phi') \in F_\tau(\xx, \alpha)} \left\{ u + \sum_{\xx' \in \mathcal{X}} V(\xx',\phi'(\xx'))\Pr[\left. \XX'=\xx' \right| \XX=\xx] \right\}, &\text{ if } F_\tau(\xx, \alpha) \neq \emptyset\\
\overline{V}, &\text{ if } F_\tau(\xx, \alpha) = \emptyset,
\end{cases}
\end{equation}
and the set $F_\tau(\xx, \alpha)$ is given by
\begin{equation}
\begin{aligned}
F_\tau(\xx, \alpha) = \Bigg\{\left(u , \phi' \right)  \Bigg| u \in \UU, \phi'(\xx') \in \Phi_\tau(\xx') \ \forall \ \xx' \in {\mathcal{X}},&\\
& \hspace{-3.0cm}\E\left[\left. c_{\tau+t+1}(u_t)\right| \XX = \xx\right] + \sum_{\xx' \in \mathcal{X}} \phi'(\xx') \Pr[\left. \XX'=\xx' \right| \XX=\xx] \leq \alpha \Bigg\}. \label{eq:F_T}
\end{aligned}
\end{equation}
    
Next, we present a detailed description of the value function and use it to compute an optimal policy for~\eqref{eq:minPSPS}.

\subsubsection{Value Function.}
\label{ssec:value_func}
We evaluate~\eqref{eq:PhiT}$-$\eqref{eq:F_T} for the objective and constraints of~\eqref{eq:minPSPS_2} to obtain
the following value function. First, for a time horizon of $\tau=1$ where the last decision is taken at day $T$, we have
\[
V_1(\xx, \alpha) = \begin{cases}
0, &\text{ if } \E\left[\left. c_{T+1}(0)\right| \XX_T = \xx\right] \leq \alpha\\
1, &\text{ if } \E\left[\left. c_{T+1}(0)\right| \XX_T = \xx\right] > \alpha \text{ and } \E\left[\left. c_{T+1}(1)\right| \XX_T = \xx\right] \leq \alpha\\
\overline{V}, &\text{ if } {b}_{1}(\xx) \equiv \min_{u \in \left\{0, 1 \right\}}\E\left[\left. c_{T+1}(u)\right| \XX_{T} = \xx\right] > \alpha.\\
\end{cases}
\]
for all $\xx \in \mathcal{X}$ and $\alpha \in \RR$. For a time horizon $\tau \geq 2$, we let 
\begin{equation}
b_\tau(\xx) = \min_{u_1 \in \left\{0, 1 \right\}}\E\left[\left. c_{T-\tau +2}(u_1)\right| \XX_{T-\tau+1} = \xx\right] + \sum_{t=2}^\tau \min_{u_t} \left\{\sum_{\xx' \in \XX} \E\left[\left. c_{T-t+2}(u_t)\right| \XX_{T-t+1} = \xx'\right] \mathbf{P}_{\xx, \xx'}^{t-1}\right\}. \label{eq:b_tau}
\end{equation}
The value function is $V_\tau(\xx, \alpha)$ if $ b_\tau(\xx) > \alpha$ or otherwise given by:
\begin{equation}
\begin{aligned}
V_\tau(\xx, \alpha) = 
\min_{\substack{u, \phi'(\xx')}} \;\; &u + \sum_{\xx' \in \mathcal{X}} \Pr[\left. \XX_{T-\tau+2}'=\xx' \right| \XX_{T-\tau+1}=\xx] V_{\tau-1}\left(\xx', \phi'(\xx') \right)\\
\text{s.t.} \quad & u \in \left\{0,1 \right\} \\
&\phi'(\xx') \in \Phi_{\tau-1}(\xx')\\
&\E\left[\left. c_{T-\tau+2}(u)\right| \XX_{T-\tau+1} = \xx\right]+\sum_{\xx' \in \mathcal{X}} \phi'(\xx') \Pr[\left. \XX_{T-\tau+2}'=\xx' \right| \XX_{T-\tau+1}=\xx] \leq \alpha,
\end{aligned}
\label{eq:value_function_problem}
\end{equation}
where
\begin{equation}
\Phi_\tau(\xx) = \left[b_{\tau}(\xx) , + \infty \right[. \label{eq:phi_tau}
\end{equation}

\subsubsection{Optimal Policy.}

We compute an optimal policy for~\eqref{eq:minPSPS_2} and, therefore, for~\eqref{eq:minPSPS} using the value function given in the previous section and~\citep{chen2007non}'s approach. Let $f_{t}:\left\{ \left(\xx, \alpha\right) \in \mathcal{X} \times \Phi_{N-t}\left(\xx\right)\right\} \mapsto \mathcal{U} \times \RR^{n}$ such that $f_{t}\left(\xx, \alpha\right) = \left(f_{t}^u\left(\xx, \alpha\right) , f_{t}^\phi\left(\xx, \alpha\right) \right) = \left(\overline{u}_{t}, \overline{\phi}_{t} \right) \in F_{T-t}(\xx, \alpha)$
where
\[
\left(\overline{u}_{t}, \overline{\phi}_{t} \right) = \argmin_{(u, \phi') \in F_{T-t}(\xx, \alpha)} \left\{ u + \sum_{\xx' \in \mathcal{X}} V_{T-t}(\xx',\phi'(\xx'))\Pr[\left. \XX_{t+1}'=\xx' \right| \XX_t=\xx] \right\},
\]
for $t\in\left\{1,2,\ldots, T\right\}$.
By~\citep[Theorem 4]{chen2007non}, if $V_T(\xx,\overline{\alpha}) < \overline{V}$, then optimal policy at time $t \in \left\{1,2,\ldots, T\right\}$ is
\[
u^\star_t = f^u_t \left( \xx_t , \alpha_{t}\right).
\]
The next round's constraint threshold is given by
\[
\alpha_{t+1} =  f_{t}^\phi \left(\xx_{t}, \alpha_{t} \right)(\xx_{t+1}),
\]
for $t=1,2,\ldots, T-1$ and $\alpha_{1} = \overline{\alpha}$, the constraint threshold specified in the problem for the whole time horizon.
Lastly, after evaluating~\eqref{eq:value_function_problem} for $\tau=1,2,\ldots, T$, we compute $u^\star_t$ and $\alpha_{t+1}$ using~\eqref{eq:policy_alpha}$-$\eqref{eq:policy_u} for $\XX_t = \xx$.

Finally, the optimal PSPS scheduling policy is presented in Proposition~\ref{prop:min_psps}. For completeness, the proof is given below.

\proof{Proof of Proposition~\ref{prop:min_psps}.}
By assumption, $V_T(\XX_1=\xx, \overline{\alpha}) < \overline{V}$ and the problem is feasible. By~\citep[Theorem 2]{chen2004dynamic} and~\citep[Theorem 4]{chen2007non}, the policy~\eqref{eq:policy_alpha}$-$\eqref{eq:policy_u} is an optimal policy for~\eqref{eq:dp_minPSPS}. The policy is, therefore, optimal for~\eqref{eq:minPSPS}.\hfill\Halmos 
\endproof

\section{Critical Peak Pricing}
\label{sec:cpp}
CPP is used to reduce peak demand by temporarily increasing electricity prices. Price increases must be declared a day ahead based on current observations. The maximum number of CPP events, $M$, is constrained by contracts between the loads and the operator. Consider, for example, Hydro-Québec's CPP program which is in effect during the Winter period from December 1$^\text{st}$ to March 31$^\text{st}$ (referred to as \emph{rate flex D})~\citep{HQCPP}. During this period, the nominal electricity price is reduced by $30\%$ from $6.08\cent/$kWh to $4.28\cent/$kWh for the first $40$ kWh multiplied by the number of days in the month, and by $22\%$ from $9.38\cent/$kWh to $7.36\cent/$kWh above this monthly consumption threshold. When a CPP is called, the price increases to $50\cent/$kWh from 6 am to 9 am and/or 4 pm to 8 pm. Finally, CPP can be called for a maximum of $100$ hours, i.e., between $25$ and $33$ times a year. The objective of CPP scheduling is, therefore, to identify the $M$ days that, without intervention, would have the highest demand. The load demand is correlated with several factors, e.g., weather parameters~\citep{hor2005analyzing,herter2007exploratory} like temperature, wind speed, precipitation, etc., day of the week~\citep{hahn2009electric}, and prior demand. 

In this section, we formulate a model for CPP scheduling based on weather and demand observations. The objective is to minimize peak demand costs. 
The model has a similar structure as PSPS Scenario 1's reformulation~\eqref{eq:min_mag_noE} from Section~\ref{sec:add_PSPS} and admits the same optimal policy. We note that~\eqref{eq:min_mag_noE} is not a model for PSPS scheduling. It is only used to derive an optimal, closed-form policy because the per-round formulations of problem~\eqref{eq:min_mag_noE} and PSPS Scenario 1 is shown to be asymptotically equivalent under certain conditions in Theorem~\ref{th:policy_and_original_prob}.

We use the same notation as in the previous sections and, for example, let $u_t =1 $ denote the decision to call for CPP during day $t+1$. We let $\XX_t \in \mathcal{X} \subset \RR^n$ be the vector collecting the weather and day of the week at time $t$. Similarly, we assume that $\XX_t$ is a Markov process with known transition probabilities for all states. Let $q_{t}: \RR^n \mapsto \RR$ be a function mapping the $n$ weather readings to an estimated peak demand level~\citep{hor2005analyzing,herter2007exploratory,hahn2009electric} at time $t$. We model the cost of supplying power to the grid as a quadratic function of the demand. This function, denoted $c^\text{power}_t:\mathcal{X}, \mapsto \RR^+$ includes, for example, generation, import, startup, and shutdown costs. Let $B_t \geq 0$, $C_t \geq 0$, and $D_t \geq 0$ be, respectively, the second, first and zeroth-order coefficient of the cost function. Let $y > 0$ be the load curtailed during a CPP event. We assume that $y$ is constant and known, e.g., an averaged historical value of total curtailment as estimated by the system operator~\citep{HQstats}. In future work, we will model $y$ as uncertain as well.
Lastly in this section, we let $\overline{a}_t \in \RR$ be the revenue loss due to high prices.

We formulate the CPP scheduling problem as:
\begin{equation}
\begin{aligned}
& \min_{\substack{u_t\\ t=1,2,\ldots,T}} & & \E \left[\sum_{t=1}^T  c^\text{power}_{t+1}\left(q_{t+1}\left(\XX_{t+1}\right)\right)(1-u_{t})+ c^\text{power}_{t+1}\left(q_{t+1}\left(\XX_{t+1}\right)-y\right)u_{t} + \overline{a}_{t+1} u_{t} \right]\\
& \text{subject to}
& & u_t \in \left\{0, 1 \right\}\\
& & & \sum_{t=1}^T u_t \leq M.
\end{aligned}
\label{eq:cpp_demand}
\end{equation}
Using the law of iterated expectation~\cite[Proposition 5.1]{ross2014first}, we obtain
\begin{equation}
\begin{aligned}
& \min_{\substack{u_t\\ t=1,2,\ldots,T}} & & \E \left[\sum_{t=1}^T  \E\left[\left. c^\text{power}_{t+1}\left(q_{t+1}\left(\XX_{t+1}\right)\right)(1-u_{t})+ c^\text{power}_{t+1}\left(q_{t+1}\left(\XX_{t+1}\right)-y\right)u_{t} + \overline{a}_{t+1} u_{t}\right| \XX_{t}\right] \right]\\
& \text{subject to}
& & u_t \in \left\{0, 1 \right\}\\
& & & \sum_{t=1}^T u_t \leq M,
\end{aligned}
\label{eq:cpp_demand}
\end{equation}
In other words, based on the observation vector at time $t$, $\XX_t$, the system operator wishes to select up to $M$ days for which the cost of supplying the peak demand is highest and in excess of the revenue losses minus the cost reduction induced by CPP. Recall that $d=T-t$. Let $v_d(w,k|\XX_{d+1})$ be the expected cumulative costs at round $d$ given the observations $\XX_{d+1}$ when the cumulative cost is $w$ and $k$ out of $M$ CPPs can still be called. The associated dynamic program is
\begin{align}
v_d(w,k|&\XX_{d+1}= \xx) = \sum_{\xs \in \mathcal{X}} \min \Big\{ v_{d-1}(w + \E\left[ \left.c_{d-1}^\text{power}\left(q_{d-1}\left(\XX_{d-1}\right)\right) \right| \XX_d = \xs \right],k|\XX_{d+1}), \label{eq:dp_cpp}\\
&v_{d-1}(w + \overline{a}_{d-1} + \E\left[ \left.c_{d-1}^\text{power}\left(q_{d-1}\left(\XX_{d-1}\right) -y \right) \right| \XX_d = \xs \right],k-1|\XX_{d+1})  \Big\} \Pr\left[\left.\XX_d = \xs \right| \XX_{d+1} = \xx \right],\nonumber
\end{align}
with the boundary conditions:
\begin{align*}
v_0(w,k|\XX_{1}= \xx) &= w + \sum_{\xs \in \mathcal{X}}\E\left[ \left.c_{-1}^\text{power}\left(q_{-1}\left(\XX_{-1}\right)\right) \right| \XX_{0} = \xs \right] \Pr\left[\left.\XX_0 = \xs \right| \XX_{1} = \xx \right]\\
v_d(w,0|\XX_{d+1}= \xx) &= w + \sum_{i=0}^d\sum_{\xs \in \mathcal{X}}\E\left[ \left.c_{i-1}^\text{power}\left(q_{i-1}\left(\XX_{i-1}\right)\right) \right| \XX_{i} = \xs \right]\mathbf{P}^{i+1}_{\xx, \xs},
\end{align*}
for all $w\in \RR$, $k=0,1, \ldots, M$ and $\xx \in \mathcal{X}$. The optimal CPP scheduling policy is given next in Proposition~\ref{prop:cpp}.
\begin{proposition}
Consider the CPP scheduling problem~\eqref{eq:cpp_demand}.
At day $d=T-t$, given the observations $\XX_d = \xx$, a CPP event is called for the following day, i.e., $u_d =1$, if
\[
\E\left[ \left.q_{d-1}\left(\XX_{d-1}\right)\right| \XX_d = \xx \right] > \frac{1}{2yB_{d-1}}\left(g_{d-1}\left(\left.k-1\right|\xx\right) - g_{d-1}\left(\left.k\right|\xx\right)+ \overline{a}_{d-1} - C_{d-1} y + B_{d-1}y^2\right).
\]
where
\begin{align*}
g_{d}(k|\XX_{d+1} =\xx) &= \sum_{\xs \in \mathcal{X}} \min \left\{ g_{d-1}(k|\xs) +  \E\left[ \left.c_{d-1}^\emph{\text{power}}\left(q_{d-1}\left(\XX_{d-1}\right)\right) \right| \XX_d = \xs \right],\right.\\
&\quad \left. g_{d-1}(k-1|\xs) + a_{d-1}+ \E\left[ \left.c_{d-1}^\emph{\text{power}}\left(q_{d-1}\left(\XX_{d-1}\right) - y \right) \right| \XX_d = \xs \right]\right\}  \Pr\left[\left.\XX_d = \xs \right| \XX_{d+1} = \xx \right],
\end{align*}
with the boundary conditions:
\begin{align*}
g_{0}(k|\xx) &=\sum_{\xs \in \mathcal{X}}\E\left[ \left.c_{-1}^\emph{\text{power}}\left(q_{-1}\left(\XX_{-1}\right)\right) \right| \XX_{0} = \xs \right] \Pr\left[\left.\XX_0 = \xs \right| \XX_{1} = \xx \right] \label{g_0d_2}\\
g_{d}(0|\xx) &= \sum_{i=0}^d\sum_{\xs \in \mathcal{X}}\E\left[ \left.c_{i-1}^\emph{\text{power}}\left(q_{i-1}\left(\XX_{i-1}\right)\right) \right| \XX_{i} = \xs \right]\mathbf{P}^{i+1}_{\xx, \xs},
\end{align*} 
for all $\xx \in \RR, k \geq 1$, and $d\geq 1$.
\label{prop:cpp}
\end{proposition} 

\proof{Proof of Proposition~\ref{prop:cpp}.}
We use the same proof technique as in Proposition~\ref{prop:policy_psps} with $s_{1,d}=s_{2,d}=0$ and replace $\E\left[ \left.f\left(\XX_{d-1}\right)\right| \XX_{d-1} = \xs \right]$ and $a_{d-1}$ with $\E\left[ \left.c_{d-1}^\text{power}\left(q_{d-1}\left(\XX_{d-1}\right)\right) \right| \XX_d = \xs \right]$ and $\overline{a}_{d-1}$, respectively. This implies that $u_d = 1$ if
\begin{align*}
g_{d-1}\left(\left.k\right|0,\xx\right) + w + \E\left[ \left.c_{d-1}^\text{power}\left(q_{d-1}\left(\XX_{d-1}\right)\right) \right| \XX_d = \xs \right] >& g_{d-1}\left(\left.k-1\right|1,\xx\right) + w + \overline{a}_{d-1} \\
&\qquad+ \E\left[ \left.f\left(\XX_{d-1} - y\right)\right| \XX_{d-1} = \xs \right].
\end{align*}
Letting $c_{d-1}^\text{power}\left(\mathbf{z}\right)= B_{d-1} \mathbf{z}^2 + C_{d-1} \mathbf{z} + D_{d-1}$ and solving for $\E\left[ \left.q_{d-1}\left(\XX_{d-1}\right)\right| \XX_d = \xs \right]$, we obtain
\[
\E\left[ \left.q_{d-1}\left(\XX_{d-1}\right)\right| \XX_d = \xs \right] > \frac{1}{2yB_{d-1}}\left(g_{d-1}\left(\left.k-1\right|1,\xx\right) - g_{d-1}\left(\left.k\right|0,\xx\right)+ \overline{a}_{d-1} - C_{d-1} y + B_{d-1}y^2\right),
\]
which completes the proof. \hfill \Halmos
\endproof

Proposition~\ref{prop:cpp} provides an optimal policy for scheduling CPPs. The policy establishes that above a peak demand provided by its right-hand term, the decision maker should call a CPP, see Figure~\ref{fig:Resultats cpp}. Our approach accounts for peak events that can be absorbed by the grid without raising prices, e.g., by cheap imports. Our CPP model differs from~\citep{chen2013optimal} as it includes the quadratic cost of the demand and the revenue loss from curtailing the load. Moreover, our policy is shown to be optimal, which was not done in~\citep{chen2013optimal}.

\section{Numerical Examples}
\label{sec:num}
 
We now provide numerical examples for PSPS Scenarios~1 and~2, and CPP scheduling policies. We restrict ourselves to a comparison between closed-form policies. As mentioned above, we do not consider PSPS Scenario 3 in this section because it does not possess useful analytical structure, and as a result is significantly more computationally intensive and, therefore, not readily implementable by system operators in comparison with the other scenarios. A closed-form policy may be obtainable via approximations. This is a topic for future work.
\subsection{PSPS}

We consider four types of weather observations to evaluate the risk of wildfire ignition: temperature, relative humidity, sustained wind, and wind gusts~\citep{PSPSfactsheer}. We use historical observations from 2011 to 2020 from the Sacramento International Airport weather station to model Northern California. We use data from 2011 to 2018 and from 2019 and 2020 as training and testing sets, respectively. We consider the months of June to September as the scheduling horizon and set $T=122$. The risk thresholds of~\eqref{eq:risk_thres} are set to greater than $30^{\circ}\textrm{C}$, lower than $20\%$, greater than $25\textrm{km/h}$, and greater than $40\textrm{km/h}$, for respectively the temperature, the relative humidity, and the sustained wind and wind gust speeds. These values are more risk-averse compared to PG\&E's from Section~\ref{sec:psps_planning}. We note that we also consider daily maximum temperature observations in the wildfire risk function but omit the low fuel moisture, red flag warning, and on-the-ground observation thresholds~\citep{PSPSpolicies} due to the lack of available historical data. The weather phenomenon's transition matrices are calculated using the training data according to eight discretized states.

We let $A_t=1\textrm{B\$}$ and $a_t=0.2\textrm{M\$}$ for all $t$. The numerical value of $A_t$ is set to be a fraction of PG\&E's liability for recent wildfires, which was in the tens of billions~\citep{abatzoglou2020population,rhodes2020balancing}. 
We set $N=10$ and $s_1=s_2=2\textrm{M\$}$. 
Finally, we set the cost adjustment $\lambda=40.5\textrm{M\$}$, i.e., 15\% of the total value of lost load (VoLL) times the average daily demand (ADD), where $\textrm{VoLL}=9000\textrm{\$/MWh}$~\citep{voll} and $\textrm{ADD}=30\textrm{GWh}$~\citep{sacramento_consumption} for the state of California and Sacramento County, respectively.

\begin{figure}[tp]
\centering
\subfloat[2019]{\includegraphics[width=0.75\columnwidth]{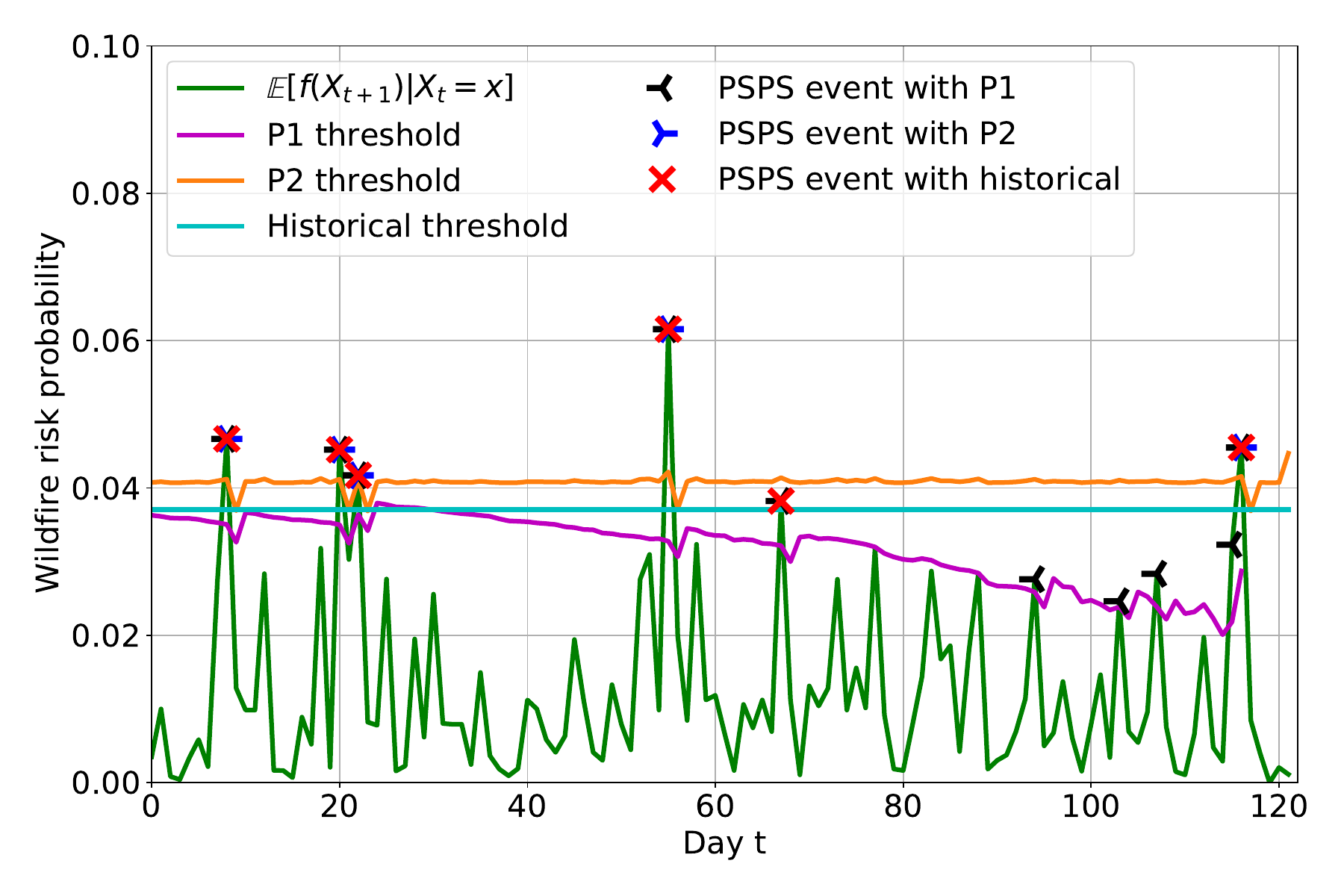}
\label{fig:PSPS 2019}
}

\subfloat[2020]{\includegraphics[width=0.75\columnwidth]{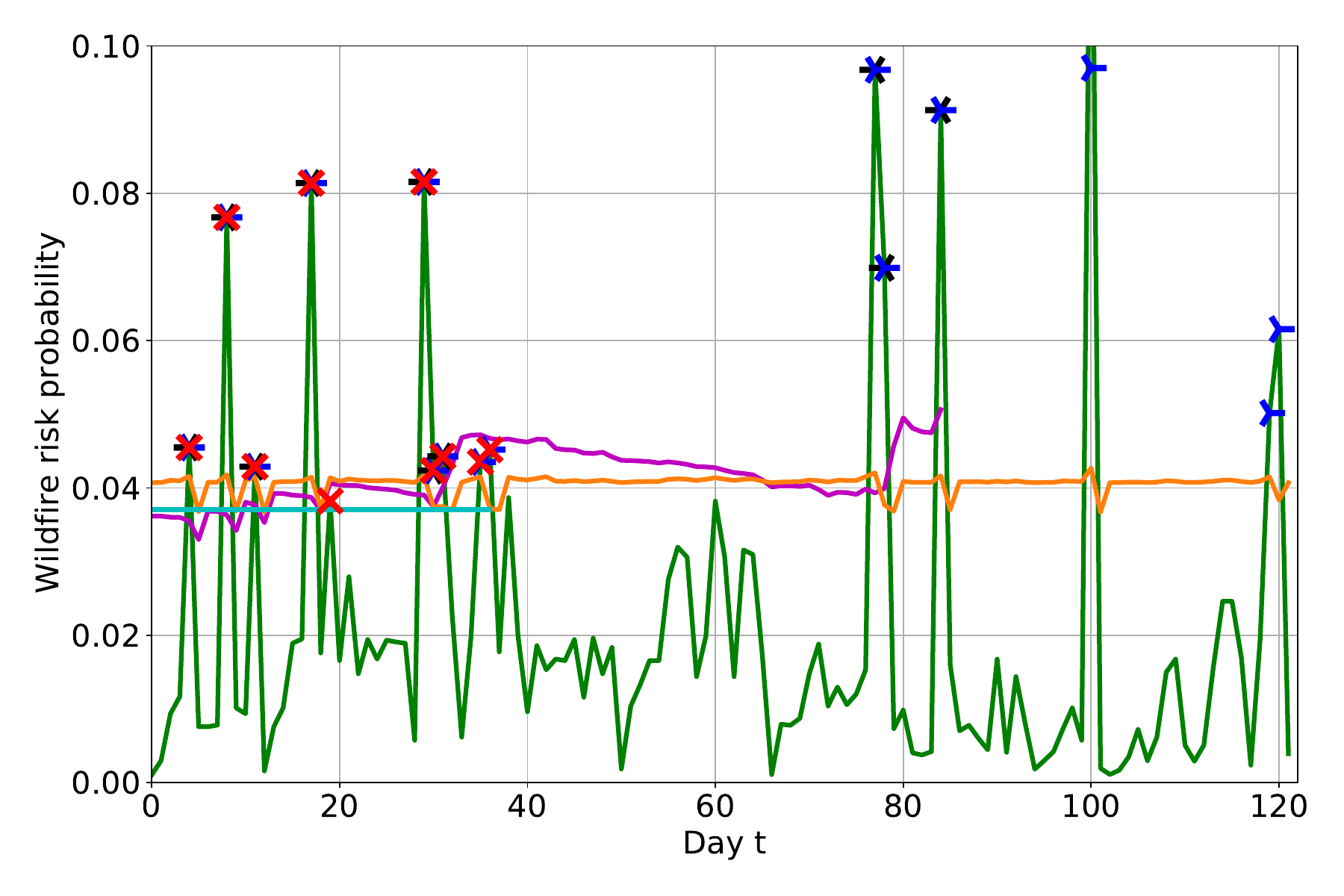}
\label{fig:PSPS 2020}
}
\caption{PSPS scheduling results for P1 and P2 \emph{($u_t = 1$ if the solid green line is above the threshold. A threshold stops when its budget is depleted.)}}
\label{fig:Resultats psps}
\end{figure}

The performances of P1 and P2, the policies for Scenarios 1 and 2, for summers 2019 and 2020 are presented in Figures~\ref{fig:PSPS 2019} and~\ref{fig:PSPS 2020}. The policies are compared to a historical policy which calls a PSPS on days with a wildfire risk probability (WRP) greater than the average $N^{\textrm{th}}$ highest WRP day for every year of the training data. The expected costs for P1, P2, and the historical policy for 2019 and 2020 are shown in Table \ref{tab:Annees test psps}. We refer to the argument of the outer expectation of, for example, Problem~\eqref{eq:min_mag_iterated} as the expected cost, viz., the cumulative conditional expected cost of each decision given the current state. Next, we apply P1, P2, and the historical policy (average 10\textsuperscript{th} highest WRP in simulated years) to 100 simulated years randomly generated using the estimated weather distribution. The average number of events and expected costs are presented in Table~\ref{tab:Annees synthetiques psps}. The policies P1 and P2 successfully select the days with the highest expected WRP as shown in Figure \ref{fig:Resultats psps}. Tables~\ref{tab:Annees test psps} and~\ref{tab:Annees synthetiques psps} show that P1 and P2 outperform the historical policy in terms of expected costs, thus leading to safer operation of the power grid.

\begin{table}[tb]
\renewcommand{\arraystretch}{1.3}
    \caption{Number of PSPS events called and expected costs for 2019 and 2020}
    \label{tab:Annees test psps}
    \centering
    \begin{tabular}{l|c|c|c|c}
        \hline
        
        \hline
        \multirow{2}{*}{\textbf{Policy}} & \multicolumn{2}{c|}{\textbf{2019}} & \multicolumn{2}{c}{\textbf{2020}}\\
        \cline{2-5}
        & Number of PSPS & Expected costs [B\$] & Number of PSPS & Expected costs [B\$]\\
        \hline
        
        \hline
        P1 & 10 & 1.140 & 10 & 1.786\\
        P2 & 5 & 1.072 & 15 & 1.660\\
        Historical & 6 & 1.240 & 10 & 1.917\\
        \hline
        
        \hline
    \end{tabular}
\end{table}

\begin{table}[tb]
\renewcommand{\arraystretch}{1.3}
    \caption{Number of PSPS events called and expected costs for 100 simulated years [avg (std)]}
    \label{tab:Annees synthetiques psps}
    \centering
    \begin{tabular}{l|c|c}
        \hline
        
        \hline
        \textbf{Policy} & \textbf{Number of PSPS} & \textbf{Expected costs [B\$]} \\
        \hline
        
        \hline
        P1 & $9.91~(0.32)$ & $1.232~(0.231)$\\
        P2 & $7.78~(3.68)$ & $1.176~(0.249)$\\
        Historical & $8.23~(2.40)$ & $1.269~(0.214)$\\
        \hline
        
        \hline
    \end{tabular}
\end{table}

\subsection{CPP}

We now present numerical results for CPP scheduling. We consider four independent weather observations: the temperatures and precipitations for both Montréal, Canada and Québec City, Canada, in addition to a week day/weekend state. We estimate $q_t$ using a linear regression~\citep{hor2005analyzing} on 2008-2018 data for the CPP period spanning December 1\textsuperscript{st} to March 31\textsuperscript{st}. We discretize $\XX$'s components and determine the number of states using leave-one-out cross-validation. We use 12 and 7 states for temperature and precipitation observations, respectively.

We set $M=25$ and consider $10^5$ clients participating to the CPP program. Each peak event is set to last 3.5 hours and we assume that participating clients reduce their power demand by $1$kW, on average~\citep{HQstats}. Accordingly, we let $a_t=3.5\textrm{h}\cdot1\textrm{kW}\cdot0.0428\textrm{\$/kWh}\cdot10^5=15\textrm{k\$}$ and $y=1\textrm{kW}\cdot10^5=100\textrm{MW}$. We set $B_t$, $C_t$, and $D_t$ to $0.00245\textrm{\$/MW\textsuperscript{2}}$, $45.5\textrm{\$/MW}$, and $800\textrm{k\$}$ for all $t$, respectively, so the generation cost may be higher than the price paid by customers on days with high demand.

\begin{figure}[tb]
\centering
\subfloat[2018-2019]{\includegraphics[width=0.75\columnwidth]{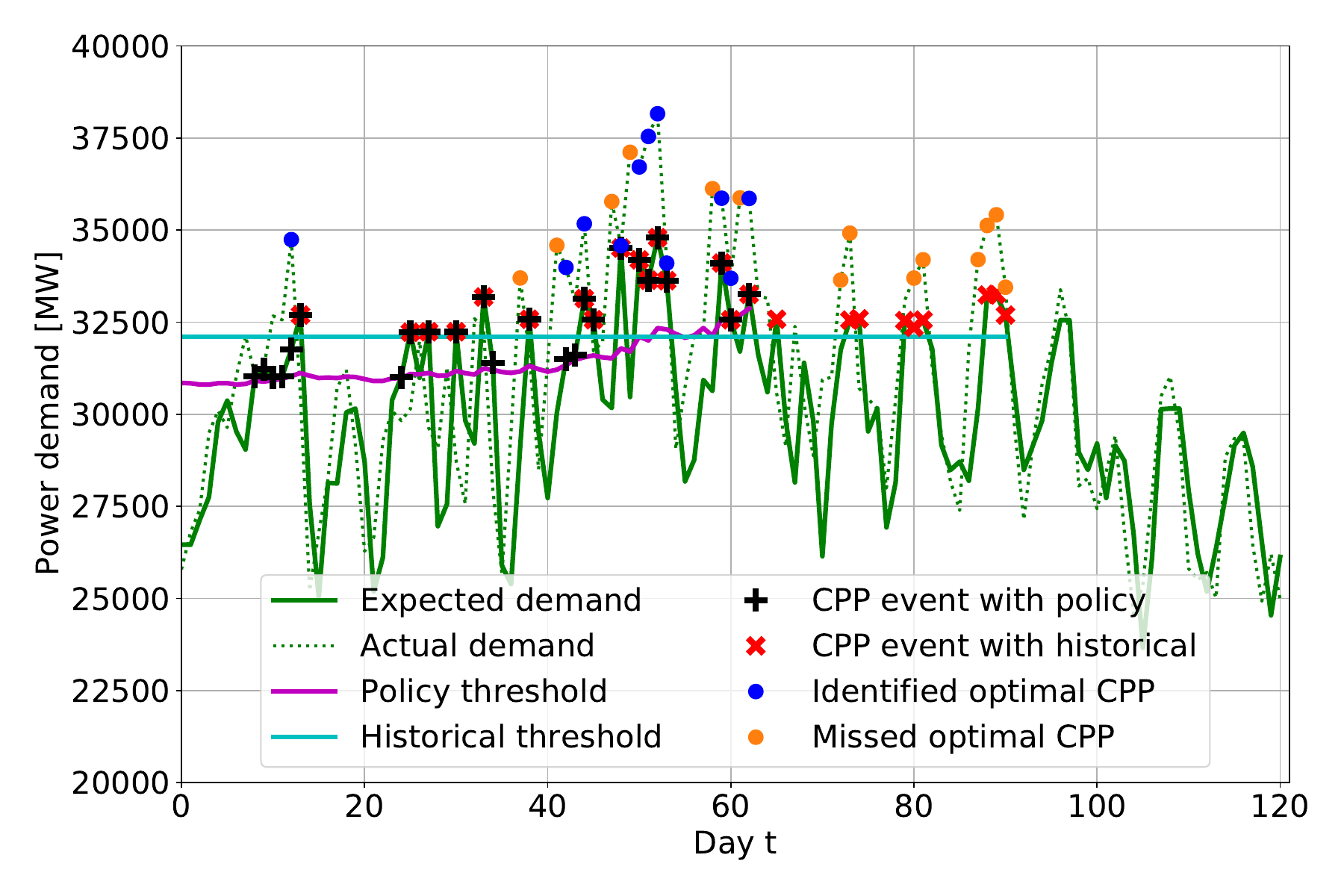}
\label{fig:CPP 2018-2019}
}

\subfloat[2019-2020]{\includegraphics[width=0.75\columnwidth]{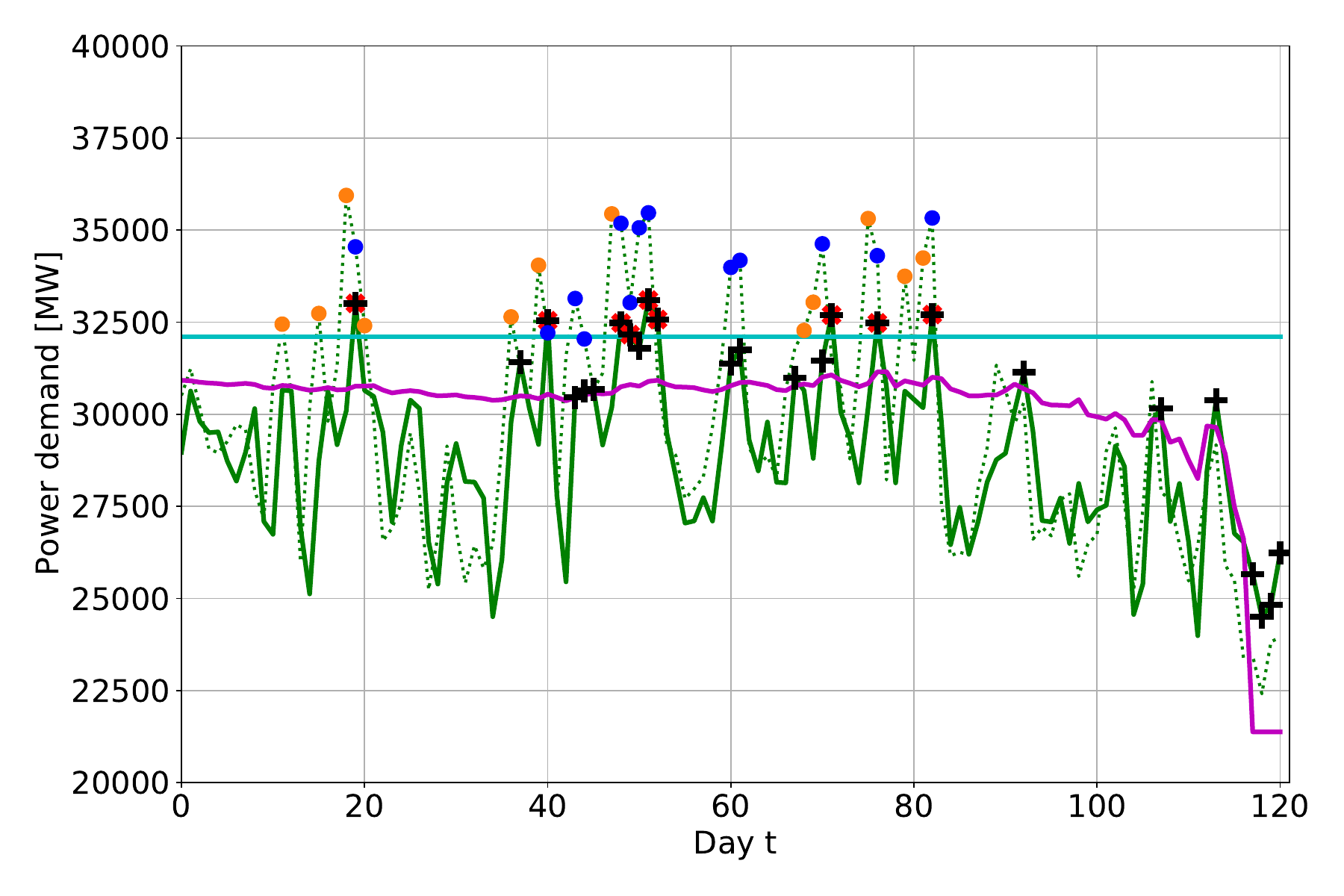}
\label{fig:CPP 2019-2020}
}
\caption{CPP scheduling results \emph{($u_t = 1$ if the solid green line is above the threshold. A threshold stops when its budget is depleted.)}}
\label{fig:Resultats cpp}
\end{figure}

The performance of our policy for 2018-2019's and 2019-2020's winters is presented in Figures~\ref{fig:CPP 2018-2019} and~\ref{fig:CPP 2019-2020}. In Figure \ref{fig:Resultats cpp}, our policy is compared to a historical policy which selects days with a demand greater than the average power demand for the $M$\textsuperscript{th} highest value for every year of the training data. The relative total cost reduction with respect to a policy selecting the days with the highest demand in hindsight is $84.1\%$ and $88.2\%$ in 2018-2019 and $78.6\%$ and $35.3\%$ in 2019-2020, respectively, for our policy and the historical policy. Figure \ref{fig:CPP 2018-2019} shows that the historical policy can outperform ours in particularly cold winters, because they differ significantly from the historical data on which our policy is based. A larger training data set and an increased number of states could address this issue. This is a topic for future work.

\begin{figure}[tb]
\centering
\includegraphics[width=0.75\columnwidth]{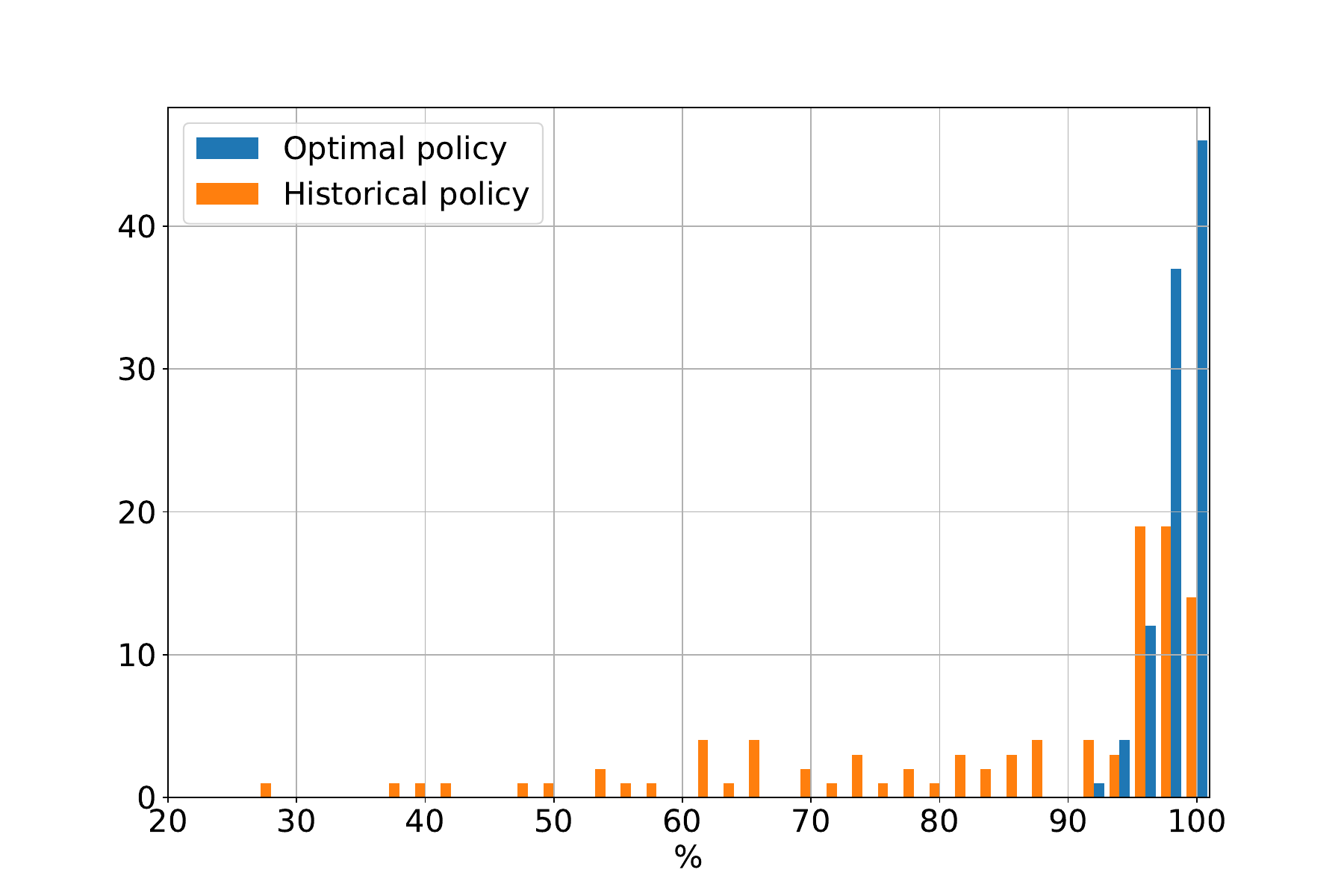}
\caption{Cost reduction distribution for 100 simulated years}
\label{fig:optimal vs simpliste}
\end{figure}

Lastly, we test our policy on $100$ simulated years. Figure~\ref{fig:optimal vs simpliste} shows the constantly high relative cost reduction with respect to the hindsight policy. Our policy leads to a $98.0\%$ average reduction with a standard deviation of $1.67\%$ and outperforms the historical policy, which leads to an $85.9\%$ average reduction with a standard deviation of $17.2\%$. The performance of our policy is thus higher and less volatile than the historical policy's.

\section{Conclusion}
\label{sec:conclu}

In this work, we formulate three dynamic programming models for PSPS scheduling to reduce power system-caused wildfires. We consider the trade-off between wildfire mitigation and the impacts of de-energizing communities by including the costs of wildfires, of PSPS operation costs, and of revenue losses for both the grid and the population. We assume that the system operator makes an initial investment to reduce wildfire risks and uses PSPSs to further decrease the risks. We consider three scenarios. First, we suppose that $N$ PSPSs are planned to reach a desired risk level and the operator must pay a penalty if the total number of PSPSs is above $N$. Second, under the same PSPS budget and penalty conditions, we assume that costs are recovered if the number of PSPS is below $N$. In the third model, the expected number of PSPSs is minimized subject to a total expected cost constraint. The first two scenarios are instances or variations of the multiple secretary problem. In each case we adapt recent results from~\cite{arlotto2019uniformly} to obtain an optimal scheduling policy for either the exact model or an asymptotically equivalent model. Lastly, we apply the first model to CPP and obtain an optimal scheduling policy. We numerically evaluate the performance of our approaches. Our simulations show that Scenario 1 and Scenario 2's policies successfully balance wildfire risk and expected costs by selecting days with the highest expected wildfire probability. P1 and P2 outperform the historical policy in test years and in simulated years. 
Lastly, the CPP scheduling policy outperforms the historical policy in simulated years, on average attaining higher savings with lower variance. However the historical policy may perform better under conditions that significantly differ from the training data, e.g., very cold winters.

Future work will focus on improving the weather model, e.g., with a larger training data set, the peak demand estimation function, and the wildfire risk probability function, e.g., with a larger number of natural phenomena, more spatial granularity, and data on actual ignitions and fire size, all of which will improve the selection of CPP and PSPS days.

\section*{Acknowledgements}
This work was funded in part by the Natural Sciences and Engineering Research Council of Canada, in part by the Institute for Data Valorization (IVADO) in part by the National Science Foundation, award 1351900, in part by the Advanced Research Projects Agency-Energy, award DE-AR0001061, and in part by the University of California Office of the President Laboratory Fees Program \#LFR-20-652467.

\bibliographystyle{informs2014}

\end{document}